\documentclass[12pt]{article}

\usepackage{amssymb}
\usepackage{amsmath}
\usepackage{amsfonts}
\usepackage{amsthm}
\usepackage{epsfig}

\theoremstyle{theorem}
    \newtheorem{theorem}{Theorem}
    \newtheorem{lemma}[theorem]{Lemma}
    \newtheorem{proposition}[theorem]{Proposition}
    
    \newtheorem{claim}[theorem]{Claim}

\theoremstyle{definition} % For roman text in the body

\theoremstyle{remark} % For an italic header, more subtle than definition style

\def\CHI{\mathchoice%
{\raise2pt\hbox{$\chi$}}%
{\raise2pt\hbox{$\chi$}}%
{\raise1.3pt\hbox{$\scriptstyle\chi$}}%
{\raise0.8pt\hbox{$\scriptscriptstyle\chi$}}}
\def\smalloplus{\raise1pt\hbox{$\,\scriptstyle \oplus\;$}}

\def\bar{\overline}

\def\<{\langle}
\def\>{\rangle}
\newcommand{\C}{{\mathbb C}}

\newcommand{\Z}{{\mathbb Z}}

\def\eqdef{\stackrel{\scriptscriptstyle def}{=}}

\def\F{{\mathcal F}}

\def \CC{{\mathcal C}}

\def\Z{\mathbb{Z}}

\def\C{\mathbb{C}}

\def\f{{\bf f}}

\def\E{{\bf E}}
\def\PP{{\mathbb P}}

\def\eqd{\stackrel{d}{=}}
\def\eqdef{\stackrel{def}{=}}

\begin{document}

%% \title[HEADLINE TITLE]{LONG TITLE \\ WITH TWO LINES}
%% The part [] is optional and can be deleted, if the title is short

\title{Large deviations for the zero set of an analytic function with diffusing coefficients}

%% First author ...name + address + email
%% \author{F. Second_name}
%% \address{Street XY, Town, State}
%% \curraddr{...}
%% current address, e-mail and url are optional
%% \email{...@...}
%% \urladdr{...}

\author{J. Ben Hough\thanks{Research partially supported by NSF grants \#DMS-0104073 and \#DMS-0244479}}
%\address{Department of Mathematics, UC Berkeley, CA 94720-3860} 
%\email{jbhough@math.berkeley.edu}
%%\urladdr{http://www.math.berkeley.edu/~jbhough}

%% If there are more authors, then second author, third author contains
%% the same items

%% (optional) If any thanks for the financial supports, grants, ...
%\thanks{Research partially supported by NSF grants \#DMS-0104073 and \#DMS-0244479}

%acknowledge the financial support from grants yyy}

%% (optional) Keywords

%% (obligatory) AMS Classification 2000
%% The Primary classification is obligatory,
%% the Secondary classification is optional.
%% \subjclass{primary}{secondary}
%% f.e. \subjclass{35R35, 49M15, 49N50}{} or \subjclass{35R35, 49M15}{49N50}

%%\subjclass{35R35, 49M15, 49N50}{}

%% (optional) Abstract

\maketitle
\begin{abstract}
The {\it hole probability} that the zero set of the time dependent planar Gaussian analytic function 
\begin{equation} \label{fc}
f_\CC(z,t) = \sum_{n=0}^\infty a_n(t) \frac{z^n}{\sqrt{n!}},
\end{equation}
where $a_n(t)$ are $i.i.d.$ complex valued Ornstein-Uhlenbeck processes, does not intersect a disk of radius $R$ for all $t \in [0,T]$ decays like $\exp(-Te^{cR^2})$.  This result sharply differentiates the zero set of $f_\CC$ from a number of canonical evolving planar point processes.  For example, the hole probability of the perturbed lattice model $\left\{\sqrt{\pi}(m,n) + c \zeta_{m,n}: m,n \in \Z \right\}$ where $\zeta_{m,n}$ are i.i.d.~Ornstein-Uhlenbeck processes decays like $\exp(-cTR^4)$.  This stark contrast is also present in the {\it overcrowding probability} that a disk of radius $R$ contains at least $N$ zeros for all $t \in [0,T]$.
\end{abstract}

%%%%% private macros, f.e. the different environments

\section{Introduction}

In this paper we study large deviations for the zero set $Z_{f_\CC}(t)$ of the time dependent planar Gaussian analytic function (GAF)
\begin{equation}
f_\CC(z,t) = \sum_{n=0}^\infty a_n(t) \frac{z^n}{\sqrt{n!}},
\end{equation}
where $a_n(t)$ are $i.i.d.$ complex valued Ornstein-Uhlenbeck processes.  Specifically, $a_n(t) = e^{-t/2}B_n(e^t)$ where $B_n(t) = \frac{1}{\sqrt{2}}\left(B_{n,1}(t) + i B_{n,2}(t)\right)$ is a Brownian motion in $\C$.  With probability one, this process defines an analytic function in the entire plane, and at any fixed time $t$ the distribution of $Z_{f_\CC}(t)$ is translation invariant (see Sodin and Tsirelson \cite{ST1} for references).  Sodin and Tsirelson \cite{ST3} study the large deviation behavior of $Z_{f_\CC}$ at a fixed time, and show that the hole probability decays exponentially in the square of the area of the disk.  The over-crowding behavior of $Z_{f_\CC}$ at a fixed time has been studied by Krishnapur \cite{Kr} who shows that the probability a fixed disk contains $N$ zeros is asymptotic to $e^{-N^2 \log N}$.  The study of Gaussian analytic functions as dynamic processes was initiated by Peres and Vir\'ag \cite{PV}, who considered the closely related time dependent hyperbolic GAF
\begin{equation} \label{hyperGAF}
f_U(z,t) = \sum_{n=0}^\infty a_n(t) z^n,
\end{equation}
where $a_n(t)$ are defined as above.  Peres and Vir\'ag showed that $|f_U(\cdot,t)|$ can be reconstructed from its zeros, and hence the zero process determined by $f_U(\cdot,t)$ is a time homogenous Markov process.  Their proof may be easily adapted to show that $|f_\CC(\cdot,t)|$ can be reconstructed from its zero set, and hence $Z_{f_\CC}(t)$ is a time homogenous Markov process as well.

For fixed $t$, the zero set $Z_{f_\CC}(t)$ exhibits strong repulsive forces between nearby zeros, as one can see visually in figure \ref{fig1}.
\begin{figure}
\centering
\includegraphics[height=1.75in]{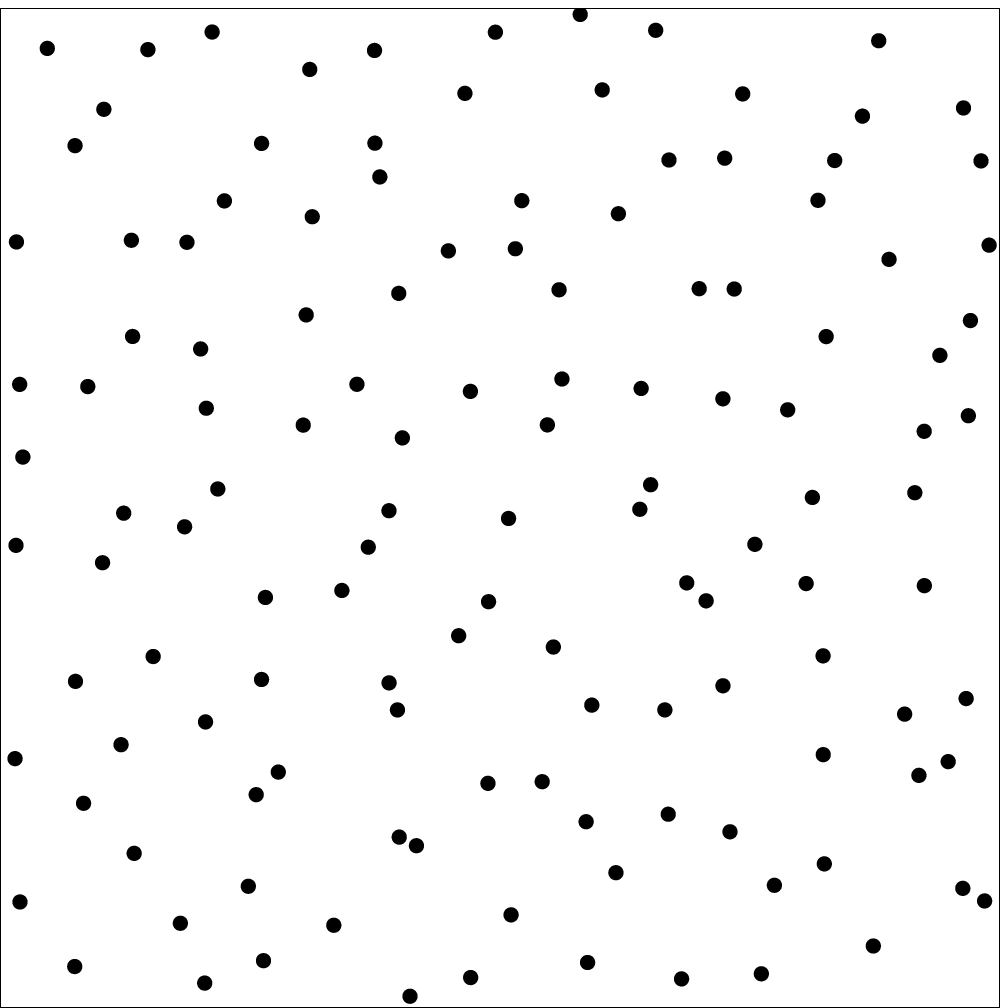}\hspace{.75in} 
\includegraphics[height=1.75in]{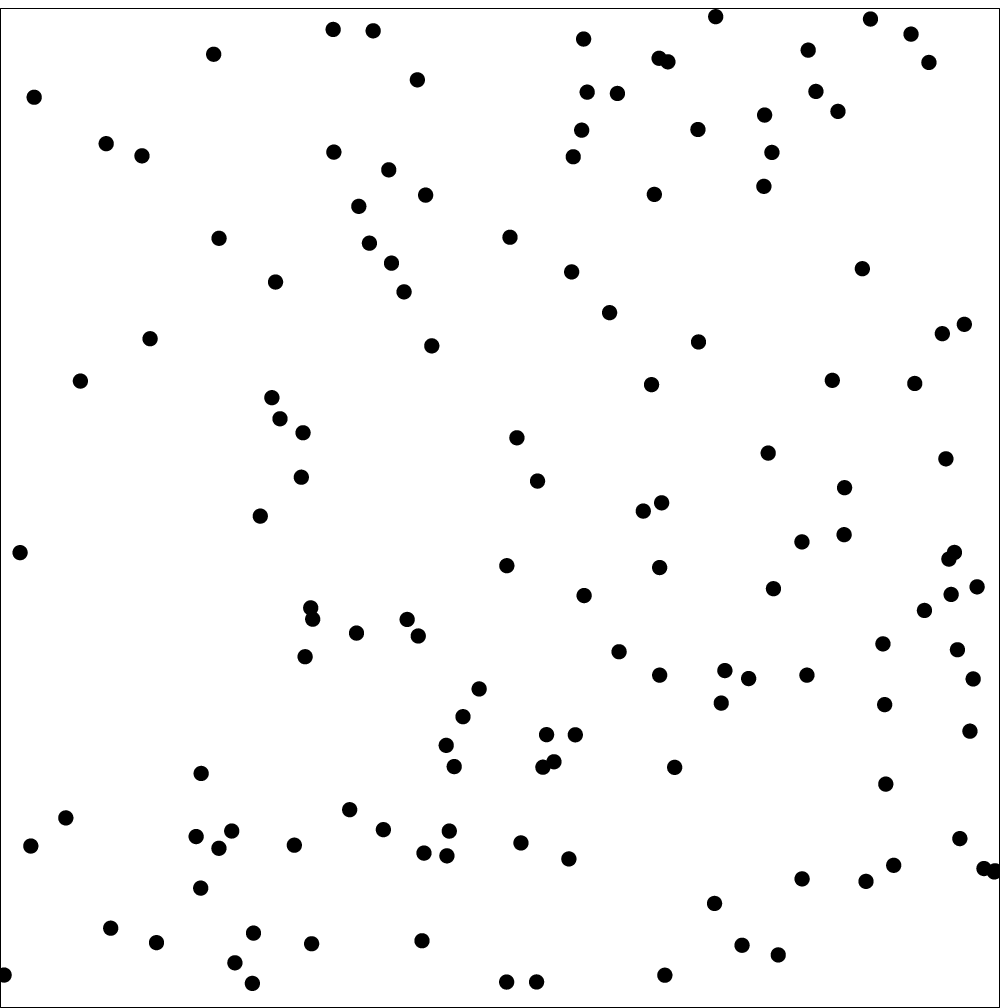}
\caption{\label{fig1} The zero set of $f_\CC(\cdot,t)$ (left) and a poisson point process with the same intensity.}
\end{figure}
To appreciate the effect of this repulsion, Sodin and Tsirelson compare $Z_{f_\CC}(0)$ to three toy models.  The first  model is a Poisson process with the same intensity, $\pi^{-1} dm$.  The second model is a {\it perturbed lattice model} consisting of the points $\left\{ \sqrt{\pi}(k + i \ell) + c \zeta_{k,\ell} : k,\ell \in  \Z \right\}$ where $\zeta_{k,\ell}$ are independent standard $\C N(0,1)$ random variables.  The third model is a {\it triangular cluster model} consisting of the points 
\begin{equation*}
\left\{ \sqrt{3\pi} (k+i\ell) + c e^{2\pi i m/3} \zeta_{k,\ell} : k,\ell \in \Z,m=0,1,2 \right\},
\end{equation*}
where $\zeta_{k,\ell}$ are as before.  They prove that the fixed time hole probability that $Z_{f_\CC}(0)$ contains no zeros in the disk of radius $r$ decays as $\exp(-c r^4)$.  This decay rate differentiates $Z_{f_{\CC}}(t)$ from the Poisson process for which the hole probability decays as $\exp(-cr^2)$, but not from the perturbed lattice and triangular cluster models.  By studying the asymptotic behavior of linear statistics, Sodin and Tsirelson are able to differentiate $Z_{f_\CC}(t)$ from the perturbed lattice model, but not the triangular cluster model.  

Each of the toy models has a natural extension to a time dependent process which preserves its distribution.  The Poisson process can be made into a time dependent process by allowing the points to evolve as independent planar Brownian motions.  The perturbed lattice and triangular cluster models may be extended to time dependent processes by replacing the coefficients $\zeta_{k,\ell}$ with independent Ornstein-Uhlenbeck processes.  We show that certain large deviation probabilities starkly contrast $Z_{f_\CC}(\cdot)$ from the other three models.  Specifically, let $H_{f_\CC}(T,R)$ denote the {\it hole probability} that $D_R = \left\{ z \in \C: | z | \leq R \right\}$ contains no zeros of $f_\CC(\cdot,t)$ for all $t \in [0,T]$.  We prove that 

\begin{theorem} \label{main} 
\begin{equation}
\limsup_{T \rightarrow \infty} \frac{1}{T} \log \left(\PP(H_{f_\CC}(T,R))\right) \leq -e^{(\frac{1}{3} - o(1))R^2}
\end{equation}
 and 
\begin{equation}
\liminf_{T \rightarrow \infty} \frac{1}{T} \log \left(\PP(H_{f_\CC}(T,R))\right) \geq -e^{(\frac{1}{2} + o(1))R^2}.
\end{equation}
\end{theorem}

Denote by $H_k(T,R)$, where $k$ equals $pp$ for Poisson process, $pl$ for perturbed lattice or $tc$ for triangular cluster, the event that the time dependent $k^{\textrm{th}}$  model has no points in $D_R$ for all $t \in [0,T]$.  We have the following results:

\begin{figure}
\centering
\includegraphics[height=1.75in]{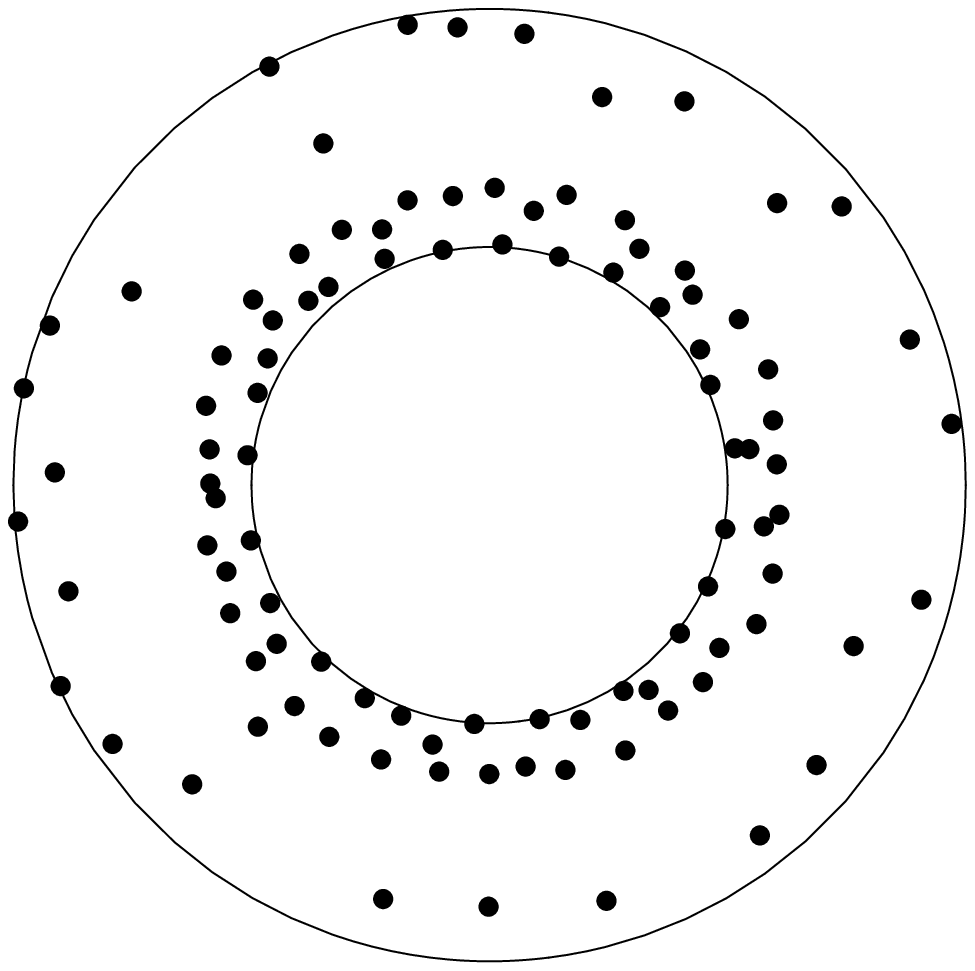}\hspace{.75in} 
\includegraphics[height=1.75in]{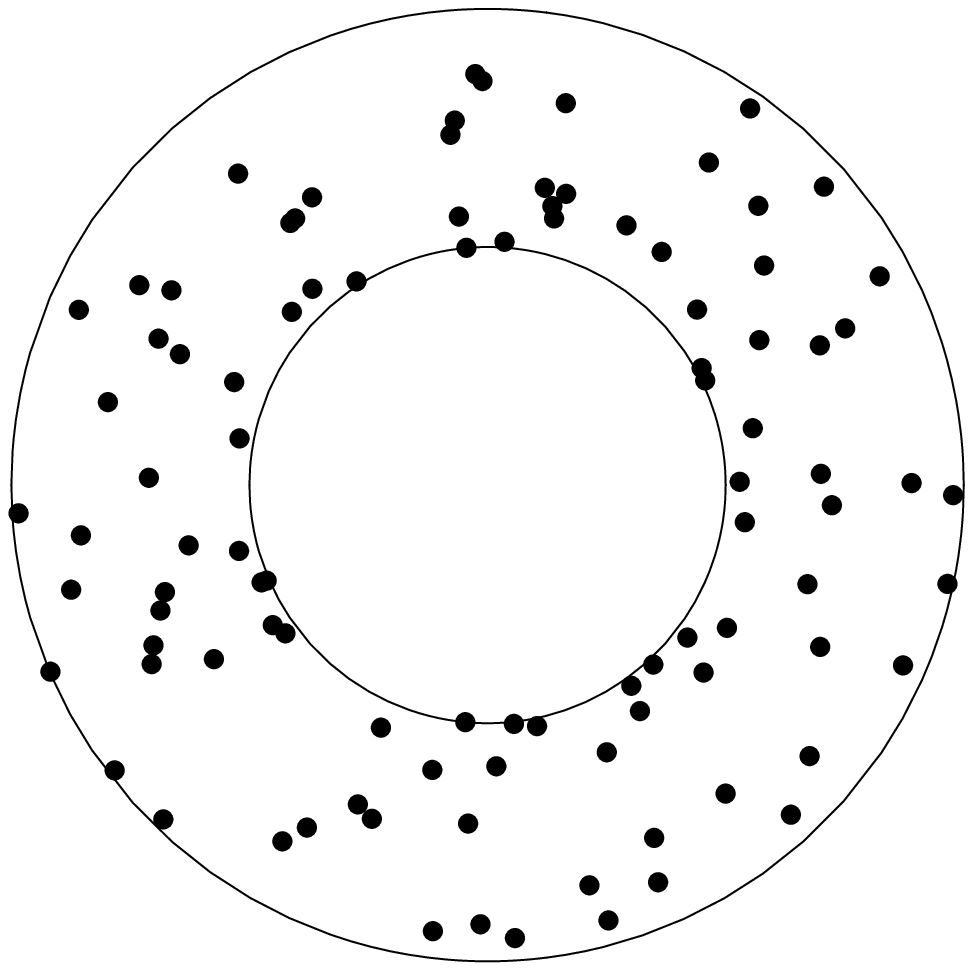}
\caption{\label{fig1} The zero set of $f_\CC(\cdot,t)$ (left) and the second toy model (right), conditioned to have a hole of radius five.}
\end{figure}

\begin{proposition} \label{toyprop1}
For the Poisson process
\begin{equation}
\lim_{T \rightarrow \infty} \frac{1}{T} \log \PP (H_{pp}(T,R)) = 0.
\end{equation}
For the perturbed lattice model and the triangular cluster model ($k = pl \textrm{ or } tc$), any $R>R_*>16$ and $T>T_*$, there exist positive constants $c_1$ and $c_2$ depending only on $T_*$ and $R_*$ so that
\begin{equation}
\limsup_{T \rightarrow \infty} \frac{1}{T} \log (\PP (H_k(T,R))) \leq -c_1 R^4 \label{toy1ub}
\end{equation}
and
\begin{equation}
 \liminf_{T \rightarrow \infty} \frac{1}{T} \log (\PP (H_k(T,R))) \geq -c_2 R^4. \label{toy1lb}
 \end{equation}
\end{proposition}

We also obtain estimates for the {\it overcrowding probability}, the event that $D_R$ contains at least $N$ zeros of $f_U(\cdot,t)$ for all $t \in [0,T]$.  Denote this event by $C_{f_\CC}(T,R,N)$.

\begin{theorem} \label{main2} For fixed $R$, we have
\begin{equation}
\limsup_{T \rightarrow \infty} \frac{1}{T} \log \left(\PP(C_{f_\CC}(T,R,N))\right) \leq -e^{(\frac{1}{6} - o(1))N \log N}
\end{equation}
 and 
\begin{equation}
\liminf_{T \rightarrow \infty} \frac{1}{T} \log \left(\PP(C_{f_\CC}(T,R,N))\right) \geq -e^{(\frac{1}{2} + o(1))N \log N}.
\end{equation}
\end{theorem}
For the three toy models, denote by $C_k(T,R,N)$ the event that the time dependent $k^{\textrm{th}}$ model ($k$ equals $pp$, $pl$ or $tc$) has at least $N$ points in $D_R$ for all $t \in [0,T]$.  We have the following results:
\begin{proposition}\label{toyprop2}
For the Poisson process 
\begin{equation}
\liminf_{T \rightarrow \infty} \frac{1}{T} \log \PP (C_{\textrm{pp}}(T,R,N)) \geq -\frac{CN}{R^2}.
\end{equation}
For the perturbed lattice model and the triangular cluster model ($k= pl \textrm{ or } tc$): 
\begin{equation}
\liminf_{T \rightarrow \infty} \frac{1}{T} \log \PP(C_k(T,R,N)) \geq -C(R)N^2.
\end{equation}
\end{proposition}

This paper is organized as follows.  In section 2, we state some well known large and small deviation estimates for Ornstein-Uhlenbeck processes, and prove the large deviation estimates for the toy models.  In section 3, we prove that $|f_\CC(z,t)|$ can be reconstructed from its zero set, and deduce that $Z_{f_\CC}(t)$ is a time homogenous Markov process.  In section 4, we prove Theorem \ref{main}, and in section 5 we prove Theorem \ref{main2}. In section \ref{OUprocess}, we prove the large deviation estimates for Ornstein-Uhlenbeck processes.

\section{Large deviations for toy models}

\subsection{Estimates for Ornstein-Uhlenbeck processes}
We begin by stating some well known large and small deviation estimates for Ornstein-Uhlenbeck processes that will be used throughout the paper.  For completeness, full proofs of these estimates are given in section \ref{OUprocess}.  These are the only estimates on the coefficients $\zeta_{k,\ell}$ necessary to derive Propositions \ref{toyprop1} and \ref{toyprop2}.  However, further properties of Ornstein-Uhlenbeck properties are used (at least superficially) to derive Theorems \ref{main} and \ref{main2}.

\begin{lemma} \label{OUsm}
Let $W(t) = e^{-t/2}B(e^t)$ where $B(t)$ is a $\C$-valued Brownian motion started from 0.  For all $R<R_*$ and $T>T_*$, there exist constants $C_1$ and $C_2$  depending only on $R_*$ and $T_*$ so that 
\begin{equation}
e^{-C_1T/R^2} \;<\; \PP( |W(t)|<R \; \forall t\in [0,T])\;<\; e^{-C_2T/R^2}.
\end{equation}
\end{lemma}

\begin{lemma}\label{OUlg}
Let $W(t) = e^{-t/2}B(e^t)$ where $B(t)$ is a $\C$-valued Brownian motion started from 0.  For all $R>R_*\geq 1$ and $T>T_*$, there exist constants $C_1$ and $C_2$ depending only on $R_*$ and $T_*$ so that 
\begin{equation}
e^{-C_1 T R^2}\; < \;\PP \left( |W(t)|>R \; \forall t\in [0,T] \right)\; < \;e^{-C_2 T R^2}.
\end{equation}
\end{lemma}

\begin{lemma}\label{OUcirc}
Let $D_\rho(x)$ denote the ball of radius $\rho$ centered at $x$.  Then for fixed $\rho$ and all $R>R_*(\rho)$ and $T>T_*(\rho)$ there exist constants $c_1$ and $c_2$ so that
\begin{equation}
e^{-c_1R^2T} \leq \PP(W(t) \in D_\rho(R) \; \forall t\in [0,T]) \leq e^{-c_2 R^2 T}.
\end{equation}
\end{lemma}

\begin{lemma} \label{OUhp} Let $W(t) = W_x(t) + i W_y(t) = e^{-t/2}B(e^t)$ where $B(t)$ is a $\C$-valued Brownian motion.  For all $R>R_*$ and $T>T_*$ there exists a constant $C$ depending only on $R_*$ and $T_*$ so that 
\begin{equation}
\PP(W_x(t) < R \; \forall t \in [0,T]) \geq \exp \left[ -T e^{-CR^2} \right].
\end{equation}
\end{lemma}

\subsection{Proof of Proposition \ref{toyprop1}}
We start with the Poisson process.  Let $\rho(\alpha,R,T)$ denote the conditional probability that no points in the Poisson process with intensity $\alpha$ intersect $D_R$ during the time interval $[0,T]$, given that no points lie in $D_R$ at time 0.  Brownian scaling and the fact that the union of two independent Poisson processes is another Poisson process gives the following:
\begin{eqnarray}
\rho(1,R,T) &=& \rho(R^2,1,T/R^2) \nonumber \\
&=& \rho(1,1, T/R^2)^{R^2}. \label{reduction}
\end{eqnarray}
Moreover, $\rho(1,1,T) = \rho(T,1/\sqrt{T},1) = \rho(1,1/\sqrt{T},1)^T$.  So it suffices to bound $\rho(1,1/\sqrt{T},1)$.

For a complex valued Brownian motion $B(t) = \frac{1}{\sqrt{2}}(B_1(t) + i B_2(t))$, let us denote by $\zeta(r)$ the hitting time of $\partial D_r$.  Recall that for $r_1<r_2<r_3$
\begin{equation}
\PP(\zeta(r_3)<\zeta(r_1) | \;|B(0)|  = r_2) = \frac{\log(r_2) - \log(r_1)}{\log(r_3) - \log(r_1)}.
\end{equation}
It follows that 
\begin{equation} \label{rad_1_est}
\PP(\zeta(\log T)< \zeta(1/\sqrt{T})|\; |B_0| = 1) = \frac{-\log (1/\sqrt{T})}{\log_2T - \log (1/\sqrt{T})} \geq 1 - \frac{2 \log_2(T)}{\log(T)}
\end{equation}
where we write $\log_2$ to denote two iterations of the $\log$ function.  Now compute
\begin{equation}
\PP \left(\max_{0\leq t \leq 1} |B(t)| \geq a \right) \leq 2 \PP\left(\max_{0 \leq t \leq 1} |B_1| \geq a \right) \leq 8 \PP(B_1(1)\geq a) \leq \frac{4 \sqrt{2}}{a \sqrt{\pi}} e^{-a^2/2}
\end{equation}
where we have used the reflection principle and Lemma 1.3 in \cite{MP}.  It follows that for large $T$ we have the estimate $\PP(\zeta(\log t)<1| \; |B(0)| = 1) \leq e^{-\frac{1}{2}(\log T)^2}$.  Combining this fact with equation (\ref{rad_1_est}) we find
\begin{equation} \label{rad1}
\PP(\zeta(1/\sqrt{T})>1 | \; |B(0)| = 1) \geq 1 - \frac{3 \log_2(T)}{\log T}.
\end{equation}
For $B(t)$ starting at radius $r<1$ we can compute the probability that $B(t)$ avoids $D_{1/\sqrt{T}}$ for all $t \in [0,1]$ by considering the probability that the Brownian motion hits $D_1$ prior to $D_{1/\sqrt{T}}$ and then use (\ref{rad1}).  This consideration yields
\begin{equation}\label{smallR}
\PP(\zeta(1/\sqrt{r}>1 | \; |B(0)| = r<1) \geq \left(1-\frac{2 \log 1/r}{\log T}\right)\left(1 - \frac{3 \log_2 T}{\log T} \right).
\end{equation}
Similar reasoning for $r>1$ yields the bound
\begin{equation}\label{bigR}
\PP(\zeta(1/\sqrt{T}>1 | \; |B(0)|=r>1) \geq 1 - \left(\frac{4\sqrt{2}}{(r-1)\sqrt{\pi}}e^{-(r-1)^2/2}\right) \left(\frac{3 \log_2T}{\log T}\right).
\end{equation}

Now fix $N$, and let $A_k=D_{(k+1)/N} \cap D_{k/N}^c$ for $1 \leq k\leq N$, and $\tilde{A}_k = D_{k+1}\cap D_k^c$ for $1\leq k \leq \infty$.  Write $\# A_k(0)$ to denote the number of points in $A_k$ at time 0, and similarly $\# \tilde{A}_k(0)$ and $(\# D_r)(0)$.  We compute, using equations (\ref{smallR}) and (\ref{bigR})
\begin{eqnarray*}
\rho(1,1/\sqrt{T},1) &\geq& \PP( (\# D_{1/N})(0) = 0)  \E  \left(\prod_{k=1}^N \left[ \left(1-\frac{2 \log\frac{N}{k}}{\log T}\right)\left(1-\frac{3 \log_2(T)}{\log T}\right) \right]^{\# A_k} \right)\\
& &\E \left(\prod_{k=1}^\infty \left[1-4(k-1) e^{-(k-1)^2/2}\left(\frac{3 \log_2 T}{\log T}\right) \right]^{\# \tilde{A}_k}\right).
\end{eqnarray*}
Now, if $M$ is a Poisson random variable with mean $\mu$, then $\E c^M = e^{(c-1)\mu}$.  Therefore
\begin{eqnarray*}
\rho(1,1/\sqrt{T},1) &\geq& e^{-1/N^2} \prod_{k=1}^N \left[\exp  \left( (-\frac{2 \log \frac{N}{k}}{\log T} - \frac{3 \log_2 T}{\log T} )  (\frac{1}{N})(\frac{2\pi(k+1)}{N})\right)\right] \\
&&\prod_{k=1}^\infty \left[\exp  \left( (- \frac{12(k-1)e^{-(k-1)^2/2} \log_2 T}{\log T})(2 \pi k) \right)\right]\\
&\geq& e^{-1/N^2} \exp \left[ \sum_{k=1}^N \left(( -\frac{2 \log \frac{N}{k}}{\log T} - \frac{3 \log_2 T}{\log T}) (\frac{1}{N})(\frac{2\pi(k+1)}{N})\right) \right]\\
 &&\exp \left[ \sum_{k=1}^\infty(- \frac{12(k-1)e^{-(k-1)^2/2} \log_2 T}{\log T})(2 \pi k) \right]\\
 &\geq& e^{-1/N^2} \exp \left[-\frac{C_1 \log_2 T}{\log T} \right] .
\end{eqnarray*}
Taking $N= \log T$, we have $\rho(1,1/\sqrt{T},1) \geq e^{- \frac{C_2 \log_2 T}{\log T}}$.  Combining this result with (\ref{reduction}) it follows that as $T \rightarrow \infty$
\begin{equation}
\frac{1}{T} \log \rho(1,R,T) \geq \frac{-C_2}{R^2} \frac{\log_2(T/R^2)}{\log (T/R^2)} \rightarrow 0,
\end{equation}
which proves the claim.

We now consider the perturbed lattice model.  To bound $\PP(H_{pl}(T,R))$ from above, note that if $H_{pl}(T,R)$ occurs then for each point $m+i n$ with $\max { |m|,|n|} \leq \lfloor R/4 \rfloor+1$ we have $|\xi_{m,n}(t) - (m,n)| \geq R/2$ for all $t \in [0,T]$.  Applying lemma \ref{OUlg}, we see that
\begin{equation}
\PP(H_2(T,R) \leq \left(e^{-C T R^2} \right)^{R^2/4}
\end{equation}
and (\ref{toy1ub}) follows.  For the lower bound, observe that $H_2(T,R)$ will occur if the following two conditions are satisfied for all lattice points $(m,n)$:
\begin{itemize}
\item[i.]  If $\max \left\{ |m|, |n| \right\} \leq 2R$ we have $|\xi_{m,n}(t) - (m+i n) | \geq 4R$ for all $t \in [0,T]$.
\item[ii.]  If $\max \left\{|m|,|n| \right\} > 2R$ then for all $t \in [0,T]$ the process $(m+ i n) + \xi_{m,n}(t)$ lies in the half plane $H$ which is a distance $R$ from the origin, is parallel to one of the coordinate axis and maximizes $d((m,n),H^c)$.
\end{itemize}
From lemma \ref{OUlg} we see that 
\begin{equation}
\PP(i) \geq \left(e^{- C T R^2} \right)^{4R^2}.
\end{equation}
From lemma \ref{OUhp} we estimate:
\begin{equation}
\PP(ii) \geq \prod_{k = \lfloor 2R \rfloor + 1}^\infty \exp \left[ - T e^{-C (k-R)^2} \right]^{8k} = e^{-\tilde{C} T}.
\end{equation}
Since the events $(i)$ and $(ii)$ are independent we obtain (\ref{toy1lb}) for $k=pl$.

\begin{figure}
\centering
\includegraphics[height=1in]{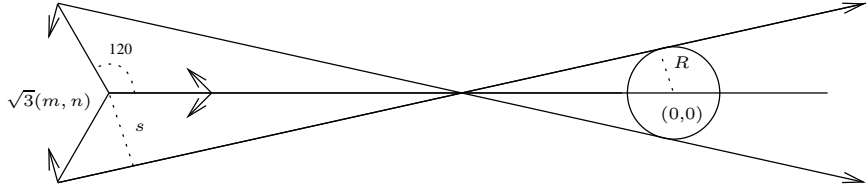}
\caption{\label{toy} Construction used for bounding the hole probability for the triangular cluster model.  By restricting one of the points in a triad to an appropriate quarter plane, we can ensure that none of the three points enters $D_R$.}
\end{figure}

The proof of the upper bound for the triangular cluster model is completely analogous.  The proof of the lower bound is also very similar.  Note that $H_{tc}(T,R)$ will occur if the following two conditions are satisfied:
\begin{itemize}
\item[i'.]  If $|\sqrt{3}(m,n)|\leq 10R$ then $|\xi_{m,n}(t) - (m,n)|\geq 20 R$ for all $t \in [0,T]$.
\item[ii'.]  If $d=|\sqrt{3}(m,n)|>10R$, then $\sqrt{3}(m,n) + c\xi_{m,n}(t)$ lies in the quarter plane $Q$  such that $Q^c$ is a distance $s = d \sin(15) - R$ from $\sqrt{3}(m,n)$ and the vertex of $Q$ lies on the segment connecting $\sqrt{3}(m,n)$ to the origin. See figure \ref{toy}.
\end{itemize}

Lemma \ref{OUhp} implies that 
\begin{equation}
\PP(i') \geq \left(e^{- C T R^2} \right)^{\tilde{C} R^2}.
\end{equation}
Moreover, since $W_x(t)$ and $W_y(t)$ are independent and the distribution of $W$ is radially symmetric we can apply lemma \ref{OUhp} to compute:
\begin{equation}
\PP(ii') \geq \prod_{k = \lfloor 10R \rfloor}^\infty \exp \left[ - 2T e^{-C_1 (k \sin(15)-R)^2} \right]^{C_2 k} = e^{-C_3 T}.
\end{equation}

\subsection{Proof of Proposition \ref{toyprop2}}
We now prove lower bounds for the over-crowding probabilities of the toy models.  
\begin{proof}[Proof of Proposition \ref{toyprop2}]
We begin with the Poisson process.  Observe that:
\begin{equation*}
\PP(B_t \in D_R \; \forall t\in[\tau,\tau+R^2] \textrm{ and } B(\tau+R^2) \in D_{R/2} | B(\tau) \in \partial D_{R/2} )
\end{equation*} 
is a constant independent of $R$, and gives a lower bound for the conditional probability
\begin{equation*}
\PP(B_t \in D_R \; \forall t\in[\tau,\tau+R^2] \textrm{ and } B(\tau+R^2) \in D_{R/2} | B(\tau) \in D_{R/2} ).
\end{equation*} 
It follows that
\begin{equation*}
\PP(C_1(N,R,T)) \geq c_1 e^{-c_2 N \lfloor T/R^2 \rfloor} 
\end{equation*}
from which the result follows.  For the perturbed lattice model, the result follows by using lemma \ref{OUcirc} to compute the probability that the points corresponding to the $N$ closest lattice points lie in $D_R$ for all $t \in [0,T]$.  The computation is analogous for the triangular cluster model, in this case we compute the probability that one point corresponding to each of the $N$ closest centers lies in $D_R$.
\end{proof}

\section{Reconstructing $|f_\CC(\cdot,t)|$ from its zero set}

This discussion closely parallels the proof given in \cite{PV} that the modulus of a hyperbolic Gaussian analytic function (\ref{hyperGAF}) can be reconstructed from its zero set.  The key lemma is the following

\begin{lemma}\label{reconstruct}
Fix $t$.  Then with probability one we have
\begin{equation} \label{keylemma}
\lim_{n \rightarrow \infty} e^{(n^6-\gamma)/2} \prod_{{z \in Z_{f_\CC}(t)} \atop{|z|<n^3}} \frac{|z|}{n^3} = |f_{\CC} (0,t)|,
\end{equation}
where $\gamma = -\int_0^\infty e^{-x} \log x \; dx$ is Euler's constant.
\end{lemma}
\begin{proof}
Assume $f_\CC(\cdot,t)$ has no zeros on $\partial D_r$, then by Jensen's formula (see \cite{Ahlf}, p. 208) we have
\begin{equation} \label{jenson_reconstruct}
\log |f_\CC(0,t)| = \frac{1}{2\pi} \int_0^{2\pi} \log |f_\CC(r e^{i \alpha})| \; d\alpha + \sum_{{z \in Z_{f_\CC}(t)}\atop{|z|<r}} \log \frac{|z|}{r}.
\end{equation}
Write $|f_\CC(re^{i \alpha},t)|^2 \eqd \sigma_r^2 Y$ where $Y \sim \exp(1)$ and $\sigma_r^2 = \textrm{Var} f_\CC(re^{i\alpha},t) = e^{r^2}$.  Then
\begin{equation*}
\E \log |f(r e^{i\alpha},t)| = \frac{1}{2}( \log (\sigma_r^2) + \E \log(Y)) = \frac{1}{2}(r^2 - \gamma).
\end{equation*}
Define $g_{r,t}(\alpha) = \log | f_\CC(re^{i \alpha},t)| + \frac{1}{2}(\gamma - r^2)$, so that the distribution of $g_{r,t}(\alpha)$ is independent of $r$, $t$ and $\alpha$, and $\E g_{r,t}(\alpha) = 0$.  Set
\begin{equation}
L_{r,t} = \frac{1}{2\pi} \int_0^{2\pi} g_{r,t}(\alpha) \; d\alpha,
\end{equation}
we prove that with probability one $L_{n^3,t} \rightarrow 0$ as $n \rightarrow \infty$.  The proof is a straightforward application of Chebyshev's inequality and the Borel-Cantelli lemma.  Compute:
\begin{eqnarray*}
\textrm{Var} L_{r,t} &=& \E \left( \frac{1}{(2\pi)^2} \int_0^{2\pi}\int_0^{2\pi} g_{r,t}(\alpha) g_{r,t}(\beta) \; d\alpha d\beta \right) \\
&=& \frac{1}{2\pi} \int_0^{2\pi} \E(g_{r,t}(\alpha) g_{r,t}(0)) \; d\alpha
\end{eqnarray*}
where we have used rotational invariance and absolute integrability to obtain the last expression. By lemma 17 in \cite{PV} we have the estimate
\begin{equation*}
\E(g_{r,t}(\alpha) g_{r,t}(0)) \leq c \frac{|\E f_\CC(r e^{i \alpha},t) \bar{f_\CC(r,t})|}{\textrm{Var}(f_\CC(r,t))} = c e^{r^2(\cos \alpha -1)}.
\end{equation*}
It follows that
\begin{eqnarray*}
\textrm{Var} L_r &\leq& \frac{1}{\pi} \int_0^\pi c e^{r^2(\cos \alpha -1)} \; d \alpha \\
&\leq& c \left[\cos^{-1}(1-1/r) + e^{-r} \right].
\end{eqnarray*}
For small $x$ we have the estimate $\cos(x) \leq 1-\frac{x^2}{2} + \frac{x^4}{24}\leq 1-\frac{11}{24}x^2$.  It follows that for $y$ close to 1, $\cos^{-1}(y) \leq \sqrt{\frac{24}{11}(1-y)}$.  Thus for large $r$ we obtain the bound
\begin{equation} \label{varest}
\textrm{Var} L_r \leq c \left[\sqrt{\frac{24}{11r}} + e^{-r} \right].
\end{equation}
Applying Chebyshev's inequality to (\ref{varest}), it follows from the Borel Cantelli lemma that $L_{n^3} \rightarrow 0$ a.s.  Thus by equation (\ref{jenson_reconstruct}):
\begin{equation}\label{keylemmaeqn}
\sum_{{z \in Z_{f_\CC}(t)} \atop{|z|<n^3}} \log \frac{|z|}{n^3} + \frac{n^6-\gamma}{2} \rightarrow \log |f_\CC(0)| \; \textrm{a.s.}
\end{equation} 
Exponentiating (\ref{keylemmaeqn}), we obtain (\ref{keylemma}).
\end{proof}

To see that $|f_\CC(\cdot,t)|$ may be reconstructed from its zero set, note that if $T(z) = z + \xi$, then by computing covariances we see that
\begin{equation*}
\tilde{f}_\CC(\cdot,t) \eqdef e^{-\bar{\xi}z - \frac{1}{2} \xi \bar{\xi}}f_\CC(T(\cdot),t) \eqd f_\CC(\cdot,t).
\end{equation*}
So, applying lemma \ref{reconstruct} to $\tilde{f}$, we can recover $|f_\CC(\xi,t)|$ with probability one.  Iterating this procedure, we can recover $|f_\CC(\cdot,t)|$ from $Z_{\f_\CC}(t)$ on a dense countable subset, and hence recover $|F_\CC(\cdot,t)|$ everywhere by continuity.

\section{Hole probability for $f_\CC$}

In this section, we compute the probability of the event $H_{f_\CC}(T,R)$ that $D_R$ contains no zeros of $f_\CC(\cdot,t)$ for all $t \in [0,T]$.   

\subsection{Proof of lower bound in Theorem \ref{main}}
We claim that a hole of radius $R$ will exist for all $t \in [0,T]$ if the following three conditions are satisfied for all such $t$:
\begin{itemize}
\item[$i.$]  $|a_0(t)| \geq 1+e^{(R^2 + \log 48R^2)/4}$ 
\item[$ii.$] $|a_k(t)| \leq e^{-(R^2 + \log 48R^2)/4}$ for $1 \leq k \leq 48R^2$
\item[$iii.$] $|a_k(t)| \leq 2^k$ for $k> 48R^2$.
\end{itemize}
A similar computation is given in \cite{ST3}.  Write $f_\CC(z,t) = a_0(t) + \psi(z,t)$ and compute
\begin{equation}
\sum_{k=1}^{48R^2} \frac{R^k |a_k|}{\sqrt{k!}} \leq \sqrt{48R^2} \sqrt{\sum_{k=1}^{48R^2} \frac{R^{2k} |a_k|^2 }{k!}} \leq \sqrt{48R^2} \max_{1\leq k \leq 48R^2}|a_k| e^{R^2/2} \leq e^{(R^2 + \log48R^2)/4}.
\end{equation}
Since we also have 
\begin{equation}
\sum_{k>48R^2} \frac{R^k}{\sqrt{k!}} 2^k \leq \sum_{k>48 R^2} \frac{2^k}{\sqrt{k!}}\left(\frac{k}{48}\right)^{k/2} < \sum_{k>48R^2}\left(\frac{k}{12} \cdot \frac{e}{k} \right)^{k/2} < \sum_{k \geq 1} 2^{-k} = \frac{1}{2}
\end{equation}
(we used the inequality $k!<\left(\frac{k}{e} \right)^k$ which follows from Sterling's formula), it follows that if ($i$), ($ii$) and ($iii$) hold then $\sup_{z \in D_R} |\psi(z,t)| \leq 1/2 + e^{(R^2 + \log48R^2)/2}$ and hence $f_\CC(z,t)$ contains no zeros in $D_R$ for $t \in [0,T]$.   Using lemmas \ref{OUsm} and \ref{OUlg} we have:
\begin{eqnarray}
	\PP(i) &\geq& \exp \left[ -T e^{(\frac{1}{2}+o(1))R^2} \right] \nonumber \\
	\PP(ii) &\geq& \exp \left[ -T e^{(\frac{1}{2}+o(1))R^2} \right]^{48R^2} = \exp \left[ -T e^{(\frac{1}{2}+o(1))R^2} \right] \nonumber\\
	\PP(iii) &\geq& \prod_{k>48R^2} e^{c_1 T/4^k} = e^{-c_2 T}.
\end{eqnarray}
Events $(i)$, $(ii)$ and $(iii)$ are independent, therefore
\begin{equation}
\PP\left(H_{f_\CC}(T,R)\right) \geq \PP(i)\PP(ii)\PP(iii) \geq \exp \left[ -T e^{(\frac{1}{2}+o(1))R^2} \right]
\end{equation}
as desired.
\subsection{Proof of upper bound in Theorem \ref{main}}
If a hole of radius $R$ exists at time $t$, Jensen's formula gives
\begin{equation}
\log |a_0(t)| = \int_{\partial D_R} \log |f_\CC(z,t)| d\mu(z),
\end{equation}
where $\mu$ is the uniform probability measure on $\partial D_R$.  So for fixed $c<1/2$ and $\tilde{c}>1/2$ one of the following three events must occur:
\begin{itemize}
\item[$A$)] $\int_{\partial D_R} \log |f_\CC(z,t)|d\mu(z) < cR^2$ and $\max_{z \in D_R} |f_\CC(z,t)| < e^{\tilde{c}R^2}$
\item[$B^1$)] $|a_0(t)| \geq e^{cR^2}$  
\item[$B^2$)] $\max_{z \in D_R} |f_\CC(z,t)| \geq e^{\tilde{c}R^2}$.
\end{itemize}
Write $A(t)$ to denote the event that $A$ occurs at time $t$, and $B^1_\ell(t)$ to denote the event that $|a_0(t)| \geq e^{(c + \ell)R^2}$ and $B^2_\ell(t)$ to denote the event that $\max_{z \in D_R} |f_\CC(z,t)| \geq e^{(\tilde{c} + \ell)R^2}$.  Let $B_\ell(t)$ denote the event that either $B^1_\ell(t)$ or $B^2_\ell(t)$ occurs.  Also, define $b(t) = \max \left\{ \ell: B_\ell(t) \textrm{ is true} \right\}$ and $\F_t = \sigma \left\{ a_k(s) \;\forall  k\geq 0 \; \textrm{and }0 \leq s\leq t\right\}$.  The method of proof is similar to the proof of the upper bound given for lemma \ref{OUlg}.  We observe the function $f_\CC(\cdot,t)$ at a sequence of times $0 = t_0<t_1<\dots< t_N$ and bound the probability that either condition $A$ or condition $B$ is satisfied at all $t_k \leq T$.  Specifically, define $t_{k+1} = t_k + \Delta t_k$, where $\Delta t_k$ is defined as follows: 
\begin{itemize}
\item[1.]  if $A(t_k)$ is true and $B_0(t_k)$ fails, $\Delta t_k = \Delta t_A(\epsilon) \eqdef e^{-(1-2c - \epsilon)R^2}$
\item[2.]  if $b(t_k) = \ell \geq 0$, then $\Delta t_k = \Delta t_B(\ell) \eqdef \left\{ \begin{array}{cc} 6R^2 & 0\leq \ell \leq 2 \\ 2(\ell + 1)R^2 & \ell>2 \end{array} \right.$.
\end{itemize}
If both $A(t_k)$ and $B_0(t_k)$ fail then we set $N=k$, i.e. the observation process is halted.  The proof relies on the following
\begin{claim} \label{claim}
For any $c>1/3$ and $0<\epsilon<1$ satisfying $\Delta t_A(\epsilon)<1$ we may choose $\tilde{c}$ and $R_*$ sufficiently large so that there exist $p_A$ and $p_B$ such that:
\begin{eqnarray}
\PP(A(t_{k+1}) \; | \; \F_{t_k}) &<& p_A \label{Aineq1} \\
\PP(B_\ell(t_{k+1}) \; | \; \F_{t_k}) &<& p_B^{\frac{\Delta t_B(\ell)}{\Delta t_A(\epsilon)}+1} \label{Bineq1}.
\end{eqnarray}
and $p_A + p_B + p_B(1-p_B)<1/2$ for all $R\geq R_*$.
\end{claim}
The proof of this claim is somewhat technical, so we shall first check that it implies the upper bound stated in Theorem \ref{main}.  

Let $p_A$ and $p_B$ be chosen as in Claim \ref{claim}, and consider the the process $\tilde{t}_n = \sum_{k=0}^{n-1} \Delta \tilde{t}_k$, where $\Delta \tilde{t}_k$ are i.i.d. and have distribution:
\begin{eqnarray}
\PP( \Delta \tilde{t}_k = 0) &=& 1 - p_A - p_B \nonumber \\
\PP(\Delta \tilde{t}_k = \Delta t_A(\epsilon)) &=& p_A + p_B(1-p_B) \nonumber \\
\PP(\Delta \tilde{t}_k = n \Delta t_A(\epsilon)) &=& p_B^n(1-p_B)\;\;\;\;\; (\textrm{for }n>1). \nonumber
\end{eqnarray}
Set $\tilde{N} = \min \left\{ k : \Delta t_k = 0 \right\}$.  Equations (\ref{Aineq1}) and (\ref{Bineq1}) imply that $\Delta \tilde{t}_k$ stochastically dominates $\Delta t_k$, so $\PP(\tilde{t}_{\tilde N} \geq T) \geq \PP(t_{N} \geq T)$.  Using the following lemma
\begin{lemma} \label{est2}
$\PP( \tilde{t}_{\tilde{N}} \geq (k+1) \Delta t_A(\epsilon) \; | \; \tilde{t}_{\tilde{N}} \geq k \Delta t_A(\epsilon)) \leq p_A + p_B + p_B(1-p_B)$.
\end{lemma}
it follows that
\begin{eqnarray*}
\PP(t_N \geq T) &\leq& \PP(\tilde{t}_{\tilde{N}} \geq \lfloor \frac{T}{\Delta t_A(\epsilon)} \rfloor \Delta t_A(\epsilon))\\
 &\leq&  (p_A + p_B + p_B(1-p_B))^{\lfloor \frac{T}{\Delta t_A(\epsilon)} \rfloor}\\
&\leq& \exp \left[ - T e^{(1-2c - \epsilon + o(1))R^2} \right],
\end{eqnarray*}
so the upper bound in Theorem \ref{main} follows by letting $c \downarrow 1/3$ and $\epsilon$ decrease to zero.  

To prove lemma \ref{est2}, we compute
\begin{eqnarray*}
\PP( \tilde{t}_{\tilde{N}} \geq (k+1) \Delta t_A(\epsilon) \; | \; \tilde{t}_{\tilde{N}} \geq k \Delta t_A(\epsilon)) &\leq& \max_{s \geq 1} \PP(\Delta \tilde{t}_n = s \Delta t_A(\epsilon) \textrm{ and } \Delta \tilde{t}_{n+1}>0 \; | \; \Delta \tilde{t}_n \geq s \Delta t_A(\epsilon)) \\
& &+ \max_{s \geq 1} \PP(\Delta \tilde{t}_n \geq (s+1) \Delta t_A(\epsilon) \; | \; \Delta \tilde{t}_n \geq s \Delta t_A(\epsilon) )  \\
&\leq& \frac{p_A + p_B - p_B^2}{p_A + p_B} (p_A + p_B) +p_B\\
&=& p_A + p_B + p_B(1-p_B).
\end{eqnarray*}

\subsection{Proof of Claim \ref{claim}}

The following lemmas allow us to bound the conditional probabilities $\PP(A(t_{k+1})|\F_{t_k})$ and $\PP(B_\ell(t_{k+1})|\F_{t_k})$.

\begin{lemma} \label{STest} For $0<\delta<1/6$, and $\tilde{z}\in \C$ with $R/2 \leq |\tilde{z}| \leq R$ and $R\geq 1$ we have:
\begin{equation}
\PP\left( \max_{z \in \tilde{z} + \delta D_R} |f_\CC(z,t + \Delta t)| \leq \sqrt{1-e^{-\Delta t}} e^{(1/2 - 3\delta) |\tilde{z}|^2 }|\;\F_t \right) \leq e^{-C_2(\delta)R^4}.
\end{equation}
\end{lemma}
\begin{proof}
Define $\|\psi\| = \sup_{z \in \tilde{z} + \delta D_R} |\psi(z)|$.  Observe that for fixed $t$ and $\Delta t$ we may write
\begin{equation}
f_\CC(\cdot,t+\Delta t) = e^{-\Delta t/2} f_\CC(\cdot,t) +  \sqrt{1-e^{-\Delta t}} q(\cdot)
\end{equation}
where $q(z) = \sum_{k=0}^\infty \frac{\alpha_k z^k}{\sqrt{k!}}$ and $\alpha_k$ are i.i.d. $\C N(0,1)$ random variables independent of $\F_t$.  Now,
\begin{equation}
\| e^{-\Delta t/2}f_\CC(\cdot,t) + \sqrt{1-e^{-\Delta t}} q(\cdot)\| + \|  e^{-\Delta t/2} f_\CC(\cdot,t) -  \sqrt{1-e^{-\Delta t}} q(\cdot) \| \geq 2\sqrt{1-e^{-\Delta t}} \|q(\cdot) \|,
\end{equation}
so
\begin{eqnarray}
\PP(\|e^{-\Delta t/2} f_\CC(\cdot,t) +  \sqrt{1-e^{-\Delta t}}q(\cdot)\|\leq k | \F_t)^2 &\leq& \PP(2\sqrt{1-e^{-\Delta t}} \|q(\cdot) \|<2k) \nonumber \\
&=& \PP(\sqrt{1-e^{-\Delta t}}\|q(\cdot)\|<k).
\end{eqnarray}
Quoting \cite[Claim 1]{ST3}, this probability is bounded above by $e^{-2C_2(\delta)R^4}$ provided that $\frac{k}{\sqrt{1-e^{-\Delta t}}}\leq e^{(1/2 - 3\delta) | \tilde{z}|^2}$.  Choosing the maximum allowable value for $k$, we obtain
\begin{equation}
\PP\left( \max_{z \in \tilde{z} + \delta D_R} |f_\CC(z,t + \Delta t)| \leq \sqrt{1-e^{-\Delta t}} e^{(1/2 - 3\delta) |\tilde{z}|^2 } \left|\; \F_t \right.\right) \leq e^{-C_2(\delta)R^4}.
\end{equation}
\end{proof}

\begin{lemma}  If $\Delta t \geq e^{-(1-2c-\epsilon)R^2}$ with $\epsilon>0$ then \label{Alemma}
\begin{equation}
\PP(A(t+\Delta t)| \; \F_t)\leq e^{-\tilde{C}(\epsilon)R^4}.
\end{equation}
\end{lemma}
\begin{proof}
The proof uses several estimates from \cite{ST3}, and we use similar notation.  Take $N = \lfloor 2\pi \delta^{-1} \rfloor$ and $z_j = \kappa R e^{2 \pi i j/N}$ where $\kappa = 1-\delta^{1/4}$ and $0<\delta<1$.  By lemma \ref{STest} we see that if $R>R_*(\delta)$ then with probability at least $1-e^{-C_3(\delta)R^4}$ we can choose $N$ points $\xi_1, \dots, \xi_N$ with $\xi_j \in z_j + \delta D_R$ such that
\begin{equation}
|f_\CC(\xi_j, t+\Delta t)| \geq \sqrt{1-e^{-\Delta t}} e^{(1/2 - 3\delta)|z_j|^2}
\end{equation}
Let $P(z, \xi)$ be the Poisson kernel for the disk $D_R$ with $|z|=R$ and $|\xi|<R$.  Define $P_j(z) = P(z, \xi_j)$.  Then if $\mu$ is the uniform probability measure on $\partial D_R$ we have
\begin{eqnarray}
(1/2-C_5 \delta^{1/4}) R^2 + \log \sqrt{1-e^{-\Delta t}} &\leq& \frac{1}{N} \sum_{j=1}^N \log |f_\CC(\xi_j,t+\Delta t)| \nonumber \\
&\leq& \int_{\partial D_R} \left( \frac{1}{N} \sum_{j=1}^N P_j(z) \right) \log | f_\CC(z,t+\Delta t) | d\mu(z) \nonumber \\
&=& \int_{\partial D_R} \left( \frac{1}{N} \sum_{j=1}^N P_j(z) - 1\right) \log | f_\CC(z, t+\Delta t)| d\mu(z) \nonumber\\
& &+ \int_{\partial D_R} \log |f_\CC(z, t+\Delta t)| d\mu(z)\label{star}
\end{eqnarray}
For the remainder of the proof we condition on the event that $\max_{z \in D_R} | f_\CC(z, t + \Delta t)| \leq e^{\tilde{c}R^2}$, since otherwise $A(t + \Delta t)$ must fail.  So, 
\begin{equation} \label{logp}
\int_{\partial D_R} \log_+ | f_\CC(z, t+\Delta t)| d\mu(z) \leq \tilde{c}R^2.
\end{equation}
Also, by applying lemma \ref{STest} with $R$ replaced by $R/2$ and $|\tilde{z}| = R/4$, we know that except on an exceptional set of measure less than $e^{-C_6(\delta)R^4}$ we have $\tilde{\xi} \in \partial D_{R/2}$ so that $|f_\CC(\tilde{\xi}, t + \Delta t)| \geq \sqrt{1-e^{-\Delta t}} \exp \left[(1/2 - 3\delta)R^2/16 \right]$.  Then
\begin{equation}
\int_{\partial D_R} \log | f_\CC(z, t+\Delta t)| P(z, \tilde{\xi}) d \mu(z) \geq \log \sqrt{1- e^{-\Delta t}} + (1/2 - 3\delta)R^2/16.
\end{equation}
An easy computation shows that $1/3 \leq P(z,\xi) \leq 3$ for $|z| = R$ and $|\xi| = R/2$, hence
\begin{eqnarray}\label{logpm}
3\int_{\partial D_R} \log_+ | f_\CC(z, t+ \Delta t)| d\mu(z) - \frac{1}{3} \int_{\partial D_R} \log_- |f_\CC(z, t+ \Delta t)| d\mu(z)\nonumber \\ \geq \int_{\partial D_R} \log |f_\CC(z, t+ \Delta t)| P(z, \tilde{\xi}) d\mu(z) \nonumber \\
\geq \log \sqrt{1-e^{-\Delta t}} + (1/2 - 3\delta) R^2/16
\end{eqnarray}
Combining (\ref{logp}) and (\ref{logpm}) we obtain
\begin{eqnarray}
\int_{\partial D_R} \log_- |f_\CC(z, t+ \Delta t)| d \mu(z) &\leq& 9\int_{\partial D_R} \log_+ | f_\CC(z, t + \Delta t)| d \mu(z) \nonumber \\
& &- 3\log \sqrt{1-e^{-\Delta t}} - 3(1/2 - 3\delta)R^2/16 \nonumber \\
&\leq& 9 \tilde{c} R^2 - 3 \log \sqrt{1-e^{- \Delta t}} \nonumber \\
\int_{\partial D_R} \left| \log |f_\CC(z, t+ \Delta t)| \right| d \mu(z) &\leq& 10 \tilde{c}R^2 - 3\log \sqrt{1-e^{-\Delta t}} \label{logabs}
\end{eqnarray}
Now from \cite[claim 2]{ST3} we know that
\begin{equation} \label{abssum}
\max_{z \in \partial D_R} \left| \frac{1}{N} \sum_{j=1}^N P_j(z) - 1 \right| < C_3 \delta^{1/2}
\end{equation}
Combining (\ref{star}), (\ref{logabs}) and (\ref{abssum}) gives that, except on an exceptional set of probability bounded by $e^{-C_2(\delta)R^4}$:
\begin{eqnarray}
\int_{\partial D_R} \log \left| f_\CC(z,t+\Delta t) \right| d\mu(z) &\geq& (1/2 - C_5 \delta^{1/4})R^2 + \log \sqrt{1-e^{-\Delta t}} \nonumber\\
& &- C_3 \delta^{1/2} \left[10 \tilde{c}R^2 - 3 \log \sqrt{1-e^{-\Delta t}} \right] \nonumber \\
&\geq& (1/2 - C_6 \delta^{1/4})R^2 + (1+C_7 \delta^{1/2}) \log \sqrt{1-e^{-\Delta t}}.\nonumber \;\;(*)
\end{eqnarray}
All we must show is that $(*)$ exceeds $cR^2$ for sufficiently small $\delta>0$ which may be chosen uniformly in $R$.  Observe that $(*)$ is increasing in $\Delta t$, so it suffices to restrict to $\Delta t = e^{-(1-2c-\epsilon)R^2}<1/2$.  Using the inequalities $1-e^{-x}>x-x^2$ and $\log(1-x)>-C_8 x$ for $0\leq x \leq 1/2$, we compute:
\begin{eqnarray}
(*) &\geq& (1/2 - C_6 \delta^{1/4})R^2 + \frac{1}{2}(1+C_7 \delta^{1/2}) (\log \Delta t + \log(1-\Delta t)) \nonumber \\
&\geq& (1/2 - C_6 \delta^{1/4})R^2 + \frac{1}{2}(1+C_7 \delta^{1/2}) (\log \Delta t- C_8 \Delta t) \nonumber \\
&\geq& (c+\frac{\epsilon}{2})R^2 + \delta^{1/4} R^2\left[-C_6 + \frac{1}{2} C_7(2c + \epsilon-1)\right] - \frac{C_8}{2}(1+C_7 \delta^{1/2}) e^{-(1-2c-\epsilon)R^2}. \label{final}
\end{eqnarray}
It is clear that for fixed $\epsilon \in (0,1-2c)$ and $R>R_*(\epsilon)$ we may choose $\delta>0$ small enough (uniformly in $R$) so that (\ref{final}) exceeds $cR^2$.
\end{proof}

%%% B lemmas

\begin{lemma} Fix $0<\epsilon<1$ small enough so that $\Delta t_A(\epsilon) < 1$.  Then for $R>1$ \label{B1lemma1}
\begin{equation}
\PP(B^1_\ell(t_{k+1}) | B_0(t_k) \textrm{ fails}) \leq  \exp \left[ - e^{(6c + 2\epsilon - 2 + o(1))R^2} \right]^{\frac{\Delta t_B(\ell)}{\Delta t_A(\epsilon)}+1}.
\end{equation}
\end{lemma}
\begin{proof}
Assuming $B_0(t_k)$ fails, write:
\begin{equation}
|a_0(t_{k+1})| \leq e^{-\Delta t/2} | a_0(t_k)| + \sqrt{1-e^{-\Delta t}} | X|
\end{equation}
where $X \sim \C N(0,1)$ is independent of $\F_{t_k}$ and $\Delta t = \Delta t_A(\epsilon)<1$.  Using the inequalities $1-x \leq e^{-x} \leq 1-x + x^2/2$ it follows that if $B^1(t_{k+1})$ is satisfied, then
\begin{equation}
\left(1-\frac{\Delta t}{2} + \frac{(\Delta t)^2}{8}\right)e^{cR^2} + \sqrt{\Delta t} |X| \geq e^{(c + \ell)R^2}
\end{equation}
and therefore, since $(\Delta t)^2<\Delta t$:
\begin{eqnarray}
|X| &\geq& \frac{1}{\sqrt{\Delta t}}\left( e^{\ell R^2} - 1 \right) e^{c R^2} + \frac{3\sqrt{\Delta t}}{8}  e^{cR^2} \nonumber \\
&\geq& \left(e^{\ell R^2} - 1\right)e^{(1/2 - \epsilon/2)R^2} + \frac{3}{8} e^{(2c + \epsilon/2 - 1/2)R^2} \eqdef Q_1(\ell).
\end{eqnarray}
Now, $\PP(|X| \geq Q_1(\ell)) = e^{-Q_1(\ell)^2}$, so 
\begin{equation}
\PP(B^1_\ell(t_{k+1}) | B_0(t_k) \textrm{ fails}) \leq  \left[\exp\left({\frac{-Q_1(\ell)^2 \Delta t_A(\epsilon)}{\Delta t_B(\ell) + \Delta t_A(\epsilon)}}\right)\right]^{\frac{\Delta t_B(\ell)}{\Delta t_A(\epsilon)}+1}
\end{equation}
The quantity $\exp\left({\frac{-Q_1(\ell)^2 \Delta t_A(\epsilon)}{\Delta t_B(\ell)+\Delta t_A(\epsilon)}}\right)$ is decreasing for $\ell\geq 0$, it follows that
\begin{eqnarray}
\PP(B^1_\ell(t_{k+1}) | B_0(t_k) \textrm{ fails}) &\leq&  \left[ \exp\left( - e^{(6c + 2 \epsilon - 2)R^2 - \log(64 R^2)} \right)\right]^{\frac{\Delta t_B(\ell)}{\Delta t_A(\epsilon)}+1} \nonumber \\
&=& \left[ \exp \left( - e^{(6c + 2 \epsilon - 2 + o(1))R^2} \right)\right]^{\frac{\Delta t_B(\ell)}{\Delta t_A(\epsilon)}+1}
\end{eqnarray}
as desired.
\end{proof}

\begin{lemma} \label{B1lemma2}Fix $0<\epsilon<1$ so that $\Delta t_A(\epsilon)<1$.  Then for $R>1$ we have
\begin{equation*}
\PP(B^1_\ell(t_{k+1}) \; | \; B_0(t_k)) \leq \exp \left[-e^{(4c+\epsilon - 1 + o(1))R^2} \right]^{\frac{\Delta t_B(\ell)}{\Delta t_A(\epsilon)}+1}.
\end{equation*}
\end{lemma}
\begin{proof}
As before, write
\begin{equation}
|a_0(t_{k+1})| \leq |a_0(t_k)| e^{-\Delta t_k/2} + \sqrt{1-e^{-\Delta t_k}} |X|
\end{equation}
where $X \sim \C N(0,1)$ is independent of $\F_{t_k}$.  Assuming $B_\ell(t_{k+1})$ occurs, we deduce:
\begin{equation}
e^{(c+\ell)R^2} \leq e^{(c-1)R^2} + |X|
\end{equation}
so $|X| \geq e^{(c + \ell)R^2} - e^{(c-1)R^2} \geq e^{(c+\ell)R^2 -1}$.  Now,
\begin{eqnarray*}
\PP\left( |X| \geq e^{(c + \ell)R^2 -1} \right)^ {\frac{\Delta t_A(\epsilon)}{\Delta t_B(\ell)+\Delta t_A(\epsilon)}} &\leq&
\PP\left( |X| \geq e^{(c + \ell)R^2 -1} \right)^ {\frac{\Delta t_A(\epsilon)}{2\Delta t_B(\ell)}}\\
&=& \exp \left[ - e^{(4c + 2\ell + \epsilon -1)R^2 -2 -\log(2\Delta t_B(\ell))} \right]
\end{eqnarray*}
is decreasing in $\ell$ for all $\ell \geq 0$, so we obtain:
\begin{equation}
\PP(B^1_\ell(t_{k+1}) \; | \; B_0(t_k)) \leq \exp \left[ - e^{(4c + \epsilon -1 + o(1))R^2} \right]^{\frac{\Delta t_B(\ell)}{\Delta t_A(\epsilon)}+1}.
\end{equation}
\end{proof}

\begin{lemma} \label{B2lemma1} Fix $0<\epsilon<1$ so that $\Delta t_A(\epsilon)<1$, then if $\tilde{c}>2$ there exists a constant $R_*>0$ so that for all $R>R_*$:
\begin{equation}
\PP(B^2_\ell(t_{k+1}) \; | \; B_0(t_k) \textrm{ fails}) \leq \exp \left[-e^{(\tilde{c} + 3c + 3 \epsilon/2 -2 + o(1))R^2} \right]^{\frac{\Delta t_B(\ell)}{\Delta t_A(\epsilon)}+1}
\end{equation}
\end{lemma}
\begin{proof}
Write $f_\CC(\cdot, t_{k+1}) = e^{-\Delta t_k/2}f_\CC(\cdot, t_k) + \sqrt{1-e^{-\Delta t_k}} q(\cdot)$ where $q(z) = \sum_{n=0}^\infty \frac{\alpha_k z^k}{\sqrt{k!}}$ with $\alpha_k \sim \C N(0,1)$ i.i.d. and independent of $\F_{t_k}$.  Also, define $\| \psi \| = \max_{z \in D_R} |\psi(z)|$ and observe that:
\begin{equation}
\| f_\CC(\cdot, t_{k+1}) \| \leq e^{-\Delta t_k/2} \| f_U(\cdot, t_k) \|  + \sqrt{1- e^{-\Delta t_k}} \| q(\cdot) \|.\label{decomp1}
\end{equation}
If $B^2_\ell(t_{k+1})$ is satisfied and $B_0(t_k)$ fails then $\Delta t_k = \Delta t_A(\epsilon)<1$, and we have (using the inequalities $e^{-\Delta t_k/2} < 1-\frac{3 \Delta t_k}{8}$ and $1-e^{-\Delta t_k}< \Delta t_k$):
\begin{eqnarray}
e^{(\tilde{c} + \ell)R^2} &\leq& e^{-\Delta t_k/2} e^{\tilde{c}R^2} + \sqrt{1- e^{-\Delta t_k}} \| q(\cdot) \| \nonumber\\
&\leq& \left(1 - \frac{3 \Delta t_k}{8} \right) e^{\tilde{c}R^2} + \sqrt{\Delta t_k} \|q(\cdot)\| \nonumber \\
\| q(\cdot)\| &\geq& \frac{1}{\sqrt{\Delta t_k}} \left( e^{\ell R^2} -1 \right) e^{\tilde{c}R^2} + \frac{3\sqrt{\Delta t_k}}{8}  e^{\tilde{c}R^2} \label{lem38} \\
&\geq&  \left( e^{\ell R^2} -1 \right) e^{(\tilde{c} + 1/2 - c - \epsilon/2)R^2} +\frac{3}{8} e^{(\tilde{c} + c + \epsilon/2- 1/2)R^2} \eqdef Q_2(\ell). \label{eq34}
\end{eqnarray}
From \cite[p. 4]{ST3} we have the estimate:
\begin{equation}
\PP(\|q(\cdot) \| > e^{(1/2 + \alpha)R^2}) < \exp \left( - e^{\alpha R^2} \right) \label{keyest}
\end{equation}
provided that $R \geq R_0(\alpha)$, where $R_0(\alpha)$ is decreasing in $\alpha$.  Since we resticted to $\tilde{c}\geq 2$ we can fix $R_*>1$ so that (\ref{keyest}) may be used to estimate the probability of (\ref{eq34}) for all $R \geq R_*$.  For such constants, we have:
\begin{eqnarray}
 \PP\left(\| q(\cdot) \| \geq Q_2(\ell) \right) ^{\frac{\Delta t_A(\epsilon)}{\Delta t_B(\ell) + \Delta t_A(\epsilon)}} & \leq& \exp \left[-e^{\log Q_2(\ell) - R^2/2} \right]^{\frac{\Delta t_A(\epsilon)}{2\Delta t_B(\ell)}} \nonumber \\
&=& \exp \left[ -\frac{3}{8} e^{(\tilde{c} + c + \epsilon/2 - 1)R^2} + (e^{\ell R^2} - 1)e^{(\tilde{c} - c - \epsilon/2)R^2} \right]^{\frac{\Delta t_A(\epsilon)}{2\Delta t_B(\ell)}} \nonumber \\
&=& \exp \left[\frac{-\frac{3}{8} e^{(\tilde{c} + 3c + 3\epsilon/2 - 2)R^2} + (e^{\ell R^2} - 1)e^{(\tilde{c} + c +\epsilon/2 -1)R^2}}{2\Delta t_B(\ell)}\right] \label{dec}.
\end{eqnarray}
It is clear that that if $R$ is sufficiently large, (\ref{dec}) is decreasing in $\ell$ for $\ell \geq 0$, so evaluating (\ref{dec}) at $\ell = 0$ we obtain
\begin{equation}
 \PP(\| q(\cdot) \| \geq Q_2(\ell) )\leq \exp \left[-e^{(\tilde{c} + 3c + 3\epsilon/2 - 2 + o(1))R^2} \right]^ {\frac{\Delta t_B(\ell)}{\Delta t_A(\epsilon)}+1}.
\end{equation}
\end{proof}

\begin{lemma} \label{B2lemma2}For $\tilde{c}\geq 2$ and $0<\epsilon<1$ such that $\Delta t_A(\epsilon)<1$ then there exists a constant $R_*$ so that for all $R>R_*$:
\begin{equation}
\PP(B^2_\ell (t_{k+1}) \; | \; B_0(t_k)) \leq \left(\exp \left[ - e^{(\tilde{c} + 2c + \epsilon - 3/2 + o(1))R^2} \right]\right)^{\frac{\Delta t_B(\ell)}{\Delta t_A(\epsilon)}+1}.
\end{equation}
\end{lemma}
\begin{proof}
Equation (\ref{decomp1}) holds as before, but now we assume that $B_0(t_k)$ is satisfied so $\Delta t_k = \Delta t_B(b(t_k))$.  If $B^2_\ell(t_{k+1})$ holds, then
\begin{equation}
e^{(\tilde{c} + \ell)R^2} \leq e^{(\tilde{c}-1)R^2} + \|q(\cdot) \|.
\end{equation}
By requiring $R_*>1$ we deduce $\|q(\cdot) \| \geq \frac{1}{2} e^{(\tilde{c} + \ell)R^2}$.  Since we have fixed $\tilde{c} \geq 2$ we can fix a constant $R_*$ so that for all $R>R_*$ equation (\ref{keyest}) may be used to write
\begin{eqnarray}
\PP\left(\| q(\cdot) \| \geq \frac{1}{2} e^{(\tilde{c} + \ell)R^2}\right)^{\frac{\Delta t_A(\epsilon)}{\Delta t_B(\ell) + \Delta t_A(\epsilon)}} &\leq& \left(\exp \left[-\frac{1}{2} e^{(\tilde{c} +\ell - 1/2)R^2} \right]\right)^{\frac{\Delta t_A(\epsilon)}{2\Delta t_B(\ell)}} \nonumber \\
&=& \exp \left[ -\frac{1}{4 \Delta t_B(\ell)}e^{(\tilde{c} + \ell + 2c + \epsilon - 3/2)R^2} \right].\label{dec2}
\end{eqnarray}
It is easy to check that (\ref{dec2}) is decreasing in $\ell$ for all $\ell\geq 0$, so
\begin{equation}
\PP(B^2_\ell (t_{k+1}) \; | \; B_0(t_k)) \leq \left(\exp \left[ - e^{(\tilde{c} + 2c + \epsilon - 3/2 + o(1))R^2} \right]\right)^{\frac{\Delta t_B(\ell)}{\Delta t_A(\epsilon)}+1} .
\end{equation}
\end{proof}
Claim \ref{claim} now follows from lemmas \ref{Alemma} - \ref{B2lemma2}.  

\section{Over-crowding probability for $f_U$}

In this section we compute the probability of the event $C_{f_\CC}(T,R,N)$ that $D_R$ contains at least $N$ zeros of $f_U(\cdot,t)$ for all $t \in [0,T]$.  

\subsection{Proof of lower bound in Theorem \ref{main2}}
We claim that for sufficiently large $N$, the disk $D_R$ will contain at least $N$ zeros for all $t \in [0,T]$ if the following three conditions are satisfied for all such $t$:
\begin{itemize}
\item[(i)] $| a_k | < \frac{R^N}{(e^{R^2} N \cdot N!)^{1/4}}$ for all $k<N$
\item[(ii)] $|a_N| \geq 2(e^{R^2} N \cdot N!)^{1/4}$
\item[(iii)] $|a_k| < 2^{k-N}$ for all $k>N$.
\end{itemize}
To see that this is the case, observe that condition $(i)$ implies that for $t \in [0,T]$:
\begin{equation}
\sum_{k=0}^{N-1} \frac{|a_k(t)| R^k}{\sqrt{k!}} < \left( \sum_{k=0}^{N-1} |a_k|^2 \right)^{1/2} \left( \sum_{k=0}^{N-1} \frac{R^{2k}}{k!} \right)^{1/2} 
<\frac{R^N e^{R^2/2}N^{1/4}}{{N!}^{1/4}}.
\end{equation}
While $(iii)$ implies that for large $N$ and $t\in [0,T]$:
\begin{eqnarray}
\sum_{k>N} \frac{|a_k(t)| R^k}{\sqrt{k!}} &\leq& \frac{2R^{N+1}}{\sqrt{(N+1)!}} \sum_{k=0}^\infty \left( \frac{2R}{\sqrt{(N+1)}}\right)^k \nonumber \\
&=& \frac{2R^{N+1}}{\sqrt{(N+1)!}} \left( \frac{1}{1 - \frac{2R}{\sqrt{N+1}}} \right)\nonumber \\
&\leq& \frac{4R^{N+1}}{\sqrt{(N+1)!}}.
\end{eqnarray}
Thus, the claim follows by comparing the functions $\frac{a_N(t) z^N}{\sqrt{N!}}$ and $f_U(z,t)$ on $\partial D_R$ and applying Rouch\'e's theorem.  From lemmas \ref{OUsm}, \ref{OUlg} and \ref{OUhp} we compute:
\begin{equation}
\PP(i) \PP(ii) \PP(iii) \geq \left(e^{-C_R T \sqrt{N^3 N!}}\right) (e^{-C T}) = \exp\left( - T e^{(\frac{1}{2}+o(1)) N \log N} \right).
\end{equation}
Since events $(i)$, $(ii)$ and $(iii)$ are independent, this computation yields the lower bound in Theorem \ref{main2}.

\subsection{Proof of upper bound in Theorem \ref{main2}}
The style of proof is identical to that given for the upperbound of Theorem \ref{main}.  Jensen's formula states that 
\begin{equation} \label{jensen}
\log |a_0(t)| = - \sum_{{|z_k|< \rho}\atop{f_U(z_k,t) = 0}} \log \left( \frac{\rho}{|z_k|} \right) + \int_{\partial D_\rho} \log |f_U(z,t)| d \mu(z),
\end{equation}
where $\mu$ is the uniform probability measure on $\partial D_\rho$.
Evaluating (\ref{jensen}) at $\rho = \sqrt{N}$, we see that if $f_U(\cdot, t)$ has at least $N$ zeros in $D_R$ and $\sqrt{N}>R$ then:
\begin{equation}
\log |a_0(t)| \leq -N\left( \frac{1}{2} \log N - \log R \right) + \int_{\partial D_{\sqrt{N}}} \log | f_U(z,t) | d \mu(z).
\end{equation}
Thus, if $C_{f_\CC}(T,R,N)$ holds, $\gamma \in (0,1/2)$ and $N>N_*(\gamma)$ it follows that at each $t \in [0,T]$ one of the following conditions must be satisfied:
\begin{itemize}
\item[A)]  $\log |a_0(t)| \leq -\frac{\gamma}{2} N \log N$
\item[B)]  $\max_{z \in \partial D_{\sqrt{N}}} \log |f_U(z,t)| \geq \frac{\gamma}{2} N \log N$.
\end{itemize}
Write $A(t)$ to denote the event that condition $A$ is satisfied at time $t$, and $B_\ell(t)$ to denote the event that $\max_{z \in \partial D_{\sqrt{N}}} \log |f_U(z,t)| \geq (\frac{\gamma}{2} + \ell) N \log N$.  Also, define $b(t) = \max \left\{ \ell: B_\ell(t) \textrm{ is true} \right\}$.  We shall observe the function $f_U(\cdot, t)$ at a sequence of times $0=t_0<t_1< \dots < t_N$, and bound the probability that either condition $A$ or condition $B$ is satisfied at all $t_k<T$.  Define $t_{k+1} = t_k + \Delta t_k$, where $\Delta t_k$ is defined as follows:
\begin{itemize}
\item[1.]  if $A(t_k)$ is true and $B_0(t_k)$ fails then $\Delta t_k = \Delta t_A(\alpha) \eqdef e^{-(\gamma - \alpha)N \log N}$
\item[2.]  if $B_0(t_k)$ is true and $b(t_k) = \ell$, then $\Delta t_k = \Delta t_B(\ell) \eqdef \left\{ \begin{array}{cc} 6 N \log N & 0\leq \ell \leq 2 \\ 2(\ell + 1)N \log N & \ell>2 \end{array} \right.$.
\end{itemize}
If both $A(t_k)$ and $B_0(t_k)$ fail, then we set $N = k$, i.e.~the observation process is halted.  We assert the following

\begin{claim} \label{claim2}
For any $\alpha \in (\frac{2 \gamma}{3}, \gamma)$ we may choose $N_*(\gamma)$ sufficiently large so that there exists $p_A$ and $p_B$ such that:
\begin{eqnarray}
\PP(A(t_{k+1}) \; | \; \F_{t_k}) &<& p_A \label{Aineq} \\
\PP(B_\ell(t_{k+1}) \; | \; \F_{t_k}) &<& p_B^{\frac{\Delta t_B(\ell)}{\Delta t_A(\alpha)}+1} \label{Bineq}
\end{eqnarray}
and $p_A + p_B + p_B(1-p_B)<1/2$ for all $N>N_*$.
\end{claim}
It then follows exactly as before that 
\begin{equation}
\PP(t_N \geq T) \leq (p_A + p_B + p_B(1-p_B))^ {\lfloor \frac{T}{\Delta t_A(\alpha)} \rfloor} \leq \exp \left[ -T e^{(\gamma - \alpha + o(1))N \log N} \right].
\end{equation}
The upper bound in Theorem \ref{main2} now follows by letting $\gamma \uparrow 1/2$ and $\alpha \uparrow 1/3$.

Claim \ref{claim2} follows from the following three lemmas:
\begin{lemma}  For $\gamma>\alpha$ we have $\PP( A(t_{k+1}) | \F_{t_k} ) \leq 2 e^{-\alpha N \log N}$
\end{lemma}
\begin{proof}
Write $a_0(t_{k+1}) = e^{-\Delta t_{k}} a_0(t_k) + \sqrt{1-e^{-\Delta t_{k}}} X$ where $X \sim \C N(0,1)$.  The probability that $A(t_{k+1})$ holds is maximized when $a_0(t_k) = 0$, so we have  
\begin{equation}
\PP( A(t_{k+1}) | \F_{t_k}) \leq \PP \left( |X| < \frac{e^{-\frac{\gamma}{2} N \log N}}{\sqrt{1-e^{-\Delta t_k}}} \right).
\end{equation}
Since $\Delta t_k \geq \Delta t_A(\alpha)$ it follows that $e^{-\Delta t_k} < 1-\frac{\Delta t_A(\alpha)}{2}$.  Hence:
\begin{eqnarray}
\PP( A(t_{k+1}) | \F_{t_k}) &\leq& \PP \left( |X| < \sqrt{2} e^{-\frac{\alpha}{2} N \log N} \right) \nonumber \\
&=& 1 - \exp \left[-2 e^{-\alpha N \log N} \right] < 2 e^{-\alpha N \log N}.
\end{eqnarray}
\end{proof}
\begin{lemma}
If $\gamma> \alpha$ there exists a constant $N_*$ so that for all $N \geq N_*$:
\begin{equation}
\PP(B_\ell(t_{k+1}) | B_0(t_k) \textrm{ fails}) \leq \exp \left[ -e^{(\frac{3}{2} \alpha - \gamma + o(1))N \log N} \right]^{\frac{\Delta t_B(\ell)}{\Delta t_A(\alpha)}+1}
\end{equation}
\end{lemma}
\begin{proof}
Define $\| \psi \| = \max_{z \in D_{\sqrt{N}}} | \psi(z) |$.  Following the proof of lemma \ref{B2lemma1}, we write
\begin{equation}
\| f_U(\cdot, t_{k+1}) \| \leq e^{-\Delta t_k/2} \| f_U(\cdot, t_k) \| + \sqrt{1-e^{-\Delta t_k}} \|q(\cdot) \| \label{decomp}
\end{equation}
and deduce that if $B_\ell(t_{k+1})$ is satisfied but $B_0(t_k)$ fails then
\begin{eqnarray}
\| q(\cdot)\| &\geq& \frac{1}{\sqrt{\Delta t_k}} \left(e^{\ell N \log N} - 1 \right) e^{\frac{\gamma}{2} N \log N} + \frac{3\sqrt{\Delta t_k}}{8}  e^{\frac{\gamma}{2} N \log N} \nonumber \\
&\geq& \left( e^{\ell N \log N} - 1 \right) e^{(\gamma - \alpha/2) N \log N} + \frac{3}{8} e^{\frac{\alpha}{2} N \log N} \eqdef Q_2(\ell)
\end{eqnarray}
Now, using (\ref{keyest}) we obtain:
\begin{eqnarray*}
\PP(\| q(\cdot) \| \geq Q_2(\ell))^{\frac{\Delta t_A(\alpha)}{\Delta t_B(\ell) + \Delta t_A(\alpha)}} &\leq& \exp \left[ -e^{\log Q_2(\ell) - N/2} \right]^{\frac{\Delta t_A(\alpha)}{2\Delta t_B(\ell)}}  \\
&=& \exp \left[ \frac{-(e^{\ell N \log N} - 1)e^{\frac{\alpha}{2} N \log N - N/2} - \frac{3}{8}e^{(\frac{3 \alpha}{2} - \gamma)N \log N - N/2}}{2\Delta t_B(\ell)} \right]
\end{eqnarray*}
It is clear that for sufficiently large $N$, the above expression is decreasing in $\ell$ for $\ell \geq 0$, so evaluating at $\ell = 0$ we obtain:
\begin{equation}
\PP( \| q(\cdot) \| \geq Q_2(\ell) ) \leq \exp \left[ -e^{(3 \alpha/2 - \gamma + o(1))N \log N} \right]^{\frac{\Delta t_B(\ell)}{\Delta t_A(\alpha)}+1}.
\end{equation}
\end{proof}

\begin{lemma} For $N \geq 2$ we have 
\begin{equation}
\PP( B_\ell(t_{k+1}) | B_0(t_k)) \leq \exp \left[ - e^{(\alpha - \gamma/2 + o(1)) N \log N} \right]^{\frac{\Delta t_B(\ell)}{\Delta t_A(\alpha)}+1}.
\end{equation} 
\end{lemma}
\begin{proof}
Equation (\ref{decomp}) holds as before, but now we assume both $B_0(t_k)$ and $B_\ell(t_{k+1})$ are satisfied so 
\begin{equation}
e^{(\frac{\gamma}{2} + \ell) N \log N} \leq e^{(\frac{\gamma}{2} - 1) N \log N} + \|q(\cdot)\|.
\end{equation}
So for $N \geq 2$ we have $\| q(\cdot) \| \geq \frac{1}{2} e^{(\frac{\gamma}{2} + \ell)N \log N}$.  Thus,
\begin{eqnarray}
\PP( \| q(\cdot) \| \geq \frac{1}{2} e^{(\frac{\gamma}{2} + \ell) N \log N} ) ^{\frac{\Delta t_A(\alpha)}{\Delta t_B(\ell) + \Delta t_A(\alpha)}} &\leq& \left( \exp \left[- \frac{1}{2} e^{(\frac{\gamma}{2} + \ell) N \log N - N/2} \right] \right)^{\frac{\Delta t_A(\alpha)}{2\Delta t_B(\ell)}} \nonumber \\
&\leq& \exp \left[ - \frac{1}{4 \Delta t_B(\ell)} e^{(\alpha - \frac{\gamma}{2} + \ell)N \log N - N/2} \right] \nonumber \\
&\leq& \exp \left[ - e^{(\alpha - \frac{\gamma}{2} + o(1)) N \log N} \right].
\end{eqnarray}

\end{proof}

\section{Proofs of Large deviations for O.U. processes} \label{OUprocess}

In this section we give proofs of the large deviation estimates for Ornstein-Uhlenbeck processes. 

\begin{proof}[Proof of lemma \ref{OUsm}]
Take $\Delta t = R^2$, and observe that:
\begin{equation}
\PP(|W(t + \Delta t)|<R \;|\; |W(t)|<R) \leq \PP(|W(t + \Delta t)|<R \;|\; W(t)=0). \label{eq1}
\end{equation}
Also, for fixed $t$ we may write $W(t + \Delta t) = e^{-\Delta t/2} W(t) + \sqrt{1 - e^{-\Delta t}} X$, where $X \sim \C N(0,1)$ is independent of $W(t)$.  Thus,
\begin{eqnarray} \label{eq4}
\PP(|W(t + \Delta t)|<R \; | \; W(t)=0) &\leq& \PP\left(|X|< \frac{R}{\sqrt{1-e^{-\Delta t}}}\right) \\
&\leq& \PP\left( |X| < \frac{R_*}{\sqrt{1-e^{-R_*^2}}} \right)= C_*
\end{eqnarray}
since $\frac{R}{\sqrt{1-e^{-R^2}}}$ is increasing in $R$.  Thus,
\begin{equation}
\PP\left( |W(t)|<R \textrm{ for } t \in [0,T]\right) \leq C_*^{\lfloor T/{\Delta t} \rfloor} = e^{-C_2 T/R^2}
\end{equation}

For the other bound, take $\Delta t = \log(1 + R^2)$ and define
\begin{eqnarray}
Q(t) &=& \textrm{ Event that } |W(s)|<R \textrm{ for }s \in [t, t+\Delta t] \textrm{ and } |W(t+\Delta t)|<R/2 \nonumber\\
\tilde{Q} &=& \textrm{ Event that } |B(s)|<R \textrm{ for }s \in [1, 1+ R^2] \textrm{ and } |B(1+R^2)|<R/2 \nonumber.
\end{eqnarray}
Observe that
\begin{equation}
\PP(Q(t) \;|\; |W(t)|<R/2) \geq  \PP(\tilde{Q} \;|\; |B(1)|=R/2) = \tilde{C}
\end{equation}
where $\tilde{C}$ is independent of $R$ by Brownian scaling.  Hence,
\begin{equation} \label{eqn8}
\PP(|W(t)|<R \textrm{ for } t\in[0,T]) \geq \PP(|W(0)|<R/2) \tilde{C}^{\lceil T/\log(1+R^2) \rceil} \geq e^{-C_1 T/R^2}.
\end{equation}

\end{proof}

\begin{proof}[Proof of lemma \ref{OUlg}:]
For the upper bound, take $\Delta t = \log(1+1/R^2)$ and define
\begin{eqnarray}
Q(t) &=& \textrm{ Event that } |W(s)|>R \textrm{ for }s \in [t, t+\Delta t] \textrm{ and } |W(t+\Delta t)|>2R \nonumber\\
\tilde{Q} &=& \textrm{ Event that } |B(s)|>R\sqrt{1+1/R^2} \textrm{ for }s \in [1, 1+1/R^2]\textrm{ and } |B(1+1/R^2)|>2R\sqrt{1+1/R^2} \nonumber\\
\tilde{Q}^\prime &=& \textrm{ Event that } |B(s)|>R^2 + 1/2 \textrm{ for }s \in [1, 2]\textrm{ and } |B(2)|>2R^2 + 1 \nonumber.
\end{eqnarray}
Now observe that:
\begin{eqnarray}
\PP(Q(t) \; |\; |W(t)|>2R) & \geq& \PP(\tilde{Q} \; | \; |B(1)|>2R) \nonumber \\
&\geq& \PP(\tilde{Q}^\prime \; | \; |B(1)|=2R^2) \geq C_* \nonumber
\end{eqnarray}
where $C_*$ is a positive constant and the last inequality follows from Brownian scaling and the inequality $\sqrt{1+x} \leq 1 + x/2$.  It follows that
\begin{equation}
\PP(|W(t)|>R \; \forall t\in[0,T]) \geq \PP(|B(1)|>2R) C_*^{\lceil T/\Delta t \rceil} \geq e^{-C_1 R^2 T}
\end{equation}

For the lower bound, we observe $W(t)$ at times $0=t_0 < t_1 < \dots<t_N$ and then bound the probability that $|W(t_k)|\geq R$ for all $t_k\leq T$.  Specifically, if we observe $W(t)$ at time $t_k$ and $|W(t_k)| = R_k \geq R$, then we set $t_{k+1} = t_k + \Delta t_k$ where $\Delta t_k = 2 +  \lfloor 2 \log(R_k/R) + 1 \rfloor$.  If $|W(t_k)|<R$ then we halt the observation process and define $N=k$.  We are interested in bounding $\PP(t_N>T)$ from above.  Observe that if $k<N$ then we may write:
\begin{equation} \label{sd1}
|W(t_{k+1})| \leq e^{-\Delta t_k/2}|W(t_k)| + \sqrt{1-e^{-\Delta t_k}} |X_k|
\end{equation}
where $X_k \sim \C N(0,1)$ is independent of $\sigma \left\{ W(t) | t<t_k \right\}$.  Thus, $|W(t_{k+1})| \leq R/e + |X_k|$ and we see that if $r \geq R$ then
\begin{equation}\label{est}
\PP\left(|W(t_{k+1})\; | \; \geq r | W(t_i) \; \forall i\leq k\right) \; \leq \; \PP\left(|X|\geq r - R/e\right) \;\leq \;\PP\left(|X| \geq r(1-1/e)\right)\; =\; e^{-r^2(1-1/e)^2}.
\end{equation}
Let
\begin{eqnarray} \label{sd2}
C_* &=& \sup_{r \geq R} \exp \left( -\frac{r^2(1-1/e)^2}{3+ \lfloor 2\log (r/R) + 1 \rfloor} \right) \nonumber \\
&\leq& \sup_{r \geq R} \exp \left( -\frac{r^2(1-1/e)^2}{4+ 2\log (r/R)} \right)\nonumber \\
&=& \exp \left( -R^2(1-1/e)^2/4 \right)
\end{eqnarray}
so that for $n \geq 1$ we have $\PP(\Delta t_{k+1} \geq n) \leq C_*^{n+1}$ by (\ref{est}) and (\ref{sd2}).  Now let $\Delta \tilde{t}_k$ be i.i.d. geometric random variables satisfying $\PP(\Delta \tilde{t}_k = n) = C_*^n(1-C_*)$ for $n = 0,1,2, \dots$.  Define $\tilde{t}_n = \sum_{k=0}^{n-1} \Delta \tilde{t}_k$ and set $\tilde{N} = \min \left\{ k: \Delta \tilde{t}_k = 0 \right\}$.  It follows that $\Delta \tilde{t}_k$ stochastically dominates the conditional distribution of $\Delta t_k$ given $\Delta t_{i}$ for all $i< k$, so $\PP(\Delta \tilde{t}_{\tilde{N}} > T ) \geq \PP(t_N > T)$.  Using the following lemma
\begin{lemma} \label{lem1}
$\PP(\Delta \tilde{t}_{\tilde{N}}>k+1 | \Delta \tilde{t}_{\tilde{N}}>k) = C_* + C_*(1-C_*)$
\end{lemma}
it follows that 
\begin{equation}
\PP(t_N > T) \leq \PP(\Delta \tilde{t}_{\tilde{N}} > \lfloor T \rfloor) = (2C_* - C_*^2)^{\lfloor T \rfloor} \leq \ e^{- C_2 T R^2}.
\end{equation}

It remains only to prove Lemma \ref{lem1}.  Let $\phi_k$ be i.i.d. Bernoulli random variables satisfying $\PP(\phi_k = 1) = C_*$.  Define $\tau_0 =0$ and $\tau_n = \min \left\{k>\tau_{n-1}:\phi_k = 0 \right\}$.  Observe that we may construct the $\phi_k$'s so that $\Delta \tilde{t}_k = \tau_{k+1} - \tau_k-1$.  Then the event $\Delta \tilde{t}_{\tilde{N}}> n$ corresponds to the event that the $\phi_k$ process yields more than $n$ 1's before 2 consecutive 0's, and the lemma follows easily.
\end{proof}

\begin{proof}[Proof of lemma \ref{OUcirc}:]
The upperbound follows from lemma \ref{OUlg}.  We now prove the lower bound.  Set $\Delta t = \log(1 + 1/R^2)$ and let

\begin{itemize}
\item[1.]  $A(t)$ denote the event that $W(t + \Delta t) \in D_{\rho/2}(R)$ and $W(s) \in D_\rho(R)$ for all $s \in[t,t+\Delta t]$.
\item[2.]  $\tilde{A}$ denote the event that $B(1+1/R^2) \in D_{\frac{\rho}{2}\sqrt{1+1/R^2}}(R\sqrt{1+1/R^2})$ and $B(s) \in D_{\rho\sqrt{s}}(R\sqrt{s})$ for all $s \in [1,1+1/R^2]$.
\end{itemize}
Then we compute:

\begin{eqnarray*}
\PP(A(t) | W(t) \in D_{\rho/2}(R)) 
&\geq& \PP(A(0) | W(0) = R-\rho/2)\\
&\geq& \PP(\tilde{A} | B(1) = R-\rho/2). 
\end{eqnarray*}
As $R \rightarrow \infty$ we have
\begin{equation}
\PP(B(s) \in D_{\rho\sqrt{s}}(R\sqrt{s}) \; \forall s \in [1,1+1/R^2] | B(1) = R-\rho/2) \rightarrow 1
\end{equation}
so for sufficiently large $R$ we obtain
\begin{eqnarray*}
\PP(A(t) | W(t) \in D_{\rho/2}(R)) 
&\geq& \frac{1}{2} \PP(B(1+1/R^2) \in D_{\frac{\rho}{2}\sqrt{1+1/R^2}}(R\sqrt{1+1/R^2}) | B(1) = R-\rho/2)\\
&=& \frac{1}{2} \PP((R-\rho/2) + \frac{1}{R} X \in D_{\frac{\rho}{2}\sqrt{1 + 1/R^2}}(R\sqrt{1+1/R^2})
\end{eqnarray*}
where $X \sim \C N(0,1)$.  Using the approximation $\sqrt{1+1/R^2}<1 + 1/(2R^2)$ we have
\begin{eqnarray*}
\PP(A(t) | W(t) \in D_{\rho/2}(R))  &\geq& \frac{1}{2} \PP(X \in D_{\frac{R\rho}{2}\sqrt{1+1/R^2}}(\frac{R\rho+1}{2}) \\
&\rightarrow& \frac{1}{2} \PP( \textrm{Re}(X)>\frac{1}{2}).
\end{eqnarray*}
Thus, $\PP(A(t) | W(t) \in D_{\rho/2}(R))$ is bounded away from zero by a positive constant.  It follows that
\begin{equation*}
\PP(W(t) \in D_\rho(R) \; \forall t\in[0,T]) \geq e^{-c \lfloor \frac{T}{\log(1+1/R^2)} \rfloor} \geq e^{-c_1 R^2T}.
\end{equation*}
\end{proof}

\begin{proof}[Proof of lemma \ref{OUhp}:]
Observe that
\begin{eqnarray*}
\PP\left(\max_{s \in [0,1]} W_x(s) < R \textrm{ and } W_x(1)<R/2 | W_x(0)=R/2\right) &\geq& \PP\left( \sup_{s \in [0,e-1]} B_x(s) \leq \left( \sqrt{e}-1 \right) \frac{R}{2} \right) \\
&=& 1 - 2 \PP \left( B_x(e-1) \geq \left( \sqrt{e}-1 \right) \frac{R}{2}  \right) \\
&\geq& 1 - e^{-c_1 R^2} \geq \exp \left[-e^{-c_2 R^2} \right]
\end{eqnarray*}
It follows that
\begin{eqnarray*}
\PP \left( \max_{s \in [0,T]} W_x(s)<R \right) &\geq& \PP(W_x(0)<R/2) \exp \left[-e^{-c_2 R^2} \right]^{\lfloor T \rfloor} \\
 &\geq& \exp\left[ -T e^{-CR^2} \right].
\end{eqnarray*}
\end{proof}

%\section{Open problems}

% or it can be used also as \begin{acknowledgment}[Acknowledgments]
\vspace{.5 in} \noindent{\bf Acknowledgements.} 
Many thanks to Yuval Peres for suggesting the problem and numerous useful discussions and comments.  I am also grateful to Mikhail Sodin for useful discussions and comments.

%%%% references

\sc \bigskip \noindent J. Ben Hough, Department of
Mathematics, U.C.\ Berkeley, CA 94720, USA. \\{\tt
jbhough@math.berkeley.edu}, {\tt
www.math.berkeley.edu/\~{}jbhough}

\end{document}